\documentclass[manuscript]{acmart}

\usepackage{amsmath,amssymb,amsthm}

\usepackage{url}

\usepackage{algorithm}

\providecommand{\headers}[2]{}

\renewcommand{\thanks}[1]{} 

\makeatletter
\newcommand{\LoadPackageUnlessAcm}[2][]{%
  \@ifclassloaded{acmart}{%
    \PackageInfo{preamble}{Skipping package `#2' because acmart class is used}%
  }{%
    \if\relax\detokenize{#1}\relax
      \usepackage{#2}%
    \else
      \usepackage[#1]{#2}%
    \fi
  }%
}
\makeatother

\LoadPackageUnlessAcm{caption}
\LoadPackageUnlessAcm{subcaption}
\LoadPackageUnlessAcm{authblk}
\LoadPackageUnlessAcm{hyperref} 

\usepackage{amstext}
\usepackage{stmaryrd}
\usepackage{algpseudocode}

\usepackage{units}
\usepackage{cancel}
\usepackage{graphicx}

\usepackage{xcolor}
\usepackage{array}
\usepackage{verbatim}
\usepackage{booktabs}
\usepackage{calc}
\usepackage{textcomp}
\usepackage{multirow}
\usepackage{tikz,tikzscale}

\usepackage{pgfplots}
\usepackage{pgfplotstable}
\usepackage{pgfkeys}
\pgfplotsset{compat=1.18}

\definecolor{deepgreen}{RGB}{0,100,0}

\newcommand{\plotwidth}{0.49\columnwidth}
\newcommand{\plotheight}{0.4\columnwidth}

\pgfplotscreateplotcyclelist{paper markers}{
  {teal,    solid, mark=*},
  {orange,   solid, mark=triangle*},
  {blue,     solid, mark=square*},
  {red,      solid, mark=pentagon*},
  {violet,     solid, mark=otimes*},
  {brown,    solid, mark=halfcircle*},
  {black,   solid, mark=diamond*},
}

\colorlet{BarColorOne}{orange!60!white}
\colorlet{BarColorTwo}{teal!80!white}
\colorlet{BarColorThree}{red!80!black}
\colorlet{BarColorFour}{blue!60!white}
\colorlet{BarColorFive}{brown!70!black}

\pgfplotsset{legend swatch/.style={area legend, draw=none}}

\pgfplotsset{
  paperplot/.style={
      width=\plotwidth,
      height=\plotheight,
      cycle list name=paper markers,
      every axis/.append style={font=\small},
      grid=major,
      thick
    }
}

\pgfplotsset{
  longplot/.style={
      width=0.8\columnwidth,
      height=\plotheight,
      cycle list name=paper markers,
      every axis/.append style={font=\small},
      grid=major,
      thick
    }
}

\usepackage[all]{xy}

\usepackage{mathtools}
\usepackage{placeins}

\usetikzlibrary{calc}
\usetikzlibrary{shadings}

\usepgfplotslibrary{fillbetween}

\usepackage[mode=multiuser,status=draft]{fixme} 
\fxsetup{innerlayout=inline}
\fxsetup{targetlayout=colorcb}

\fxusetheme{color}
\FXRegisterAuthor{mw}{envmw}{MW}

\definecolor{apply}{rgb}{0.3,.7,0.}
\definecolor{applyface}{rgb}{.9,.95,0.}
\definecolor{invert}{rgb} {0.5,0.,1.}
\colorlet{invertface}{invert!60!red}

\usepackage[colorinlistoftodos,prependcaption,textsize=small]{todonotes}

\def\mesh{\mathbb M}

\def\R{\mathbb R}

\def\resth#1{I_{#1}^{\downarrow_{h}}}
\def\prolh#1{I_{#1}^{\uparrow_{h}}}

\def\prolp#1{I_{#1}^{\Uparrow_{p}}}
\def\restp#1{I_{#1}^{\Downarrow_{p}}}

\usepackage{colortbl} 
\newcommand{\myrowcolor}{\rowcolor[gray]{0.925}}

\usepackage{booktabs} 
\usepackage[font=footnotesize,position=b]{subcaption}
\usepackage{adjustbox}

\newcommand{\pluseq}{\stackrel{+}{\gets}}

\newcommand{\LUpNew}{local\_update}

\newcommand{\residual}{r}
\newcommand{\res}{\text{res}}

\makeatletter
\def\pgfplots@getautoplotspec into#1{%
    \begingroup
    \let#1=\pgfutil@empty
    \pgfkeysgetvalue{/pgfplots/cycle multi list/@dim}\pgfplots@cycle@dim
    \let\pgfplots@listindex=\pgfplots@numplots
    \pgfkeysgetvalue{/pgfplots/cycle list set}\pgfplots@listindex@set
    \ifx\pgfplots@listindex@set\pgfutil@empty
    \else
      \c@pgf@counta=\pgfplots@listindex
      \c@pgf@countb=\pgfplots@listindex@set
      \advance\c@pgf@countb by -\c@pgf@counta
      \globaldefs=1\relax
      \edef\setshift{%
        \noexpand\pgfkeys{
          /pgfplots/cycle list shift=\the\c@pgf@countb,
          /pgfplots/cycle list set=
        }
      }%
      \setshift%
      \globaldefs=0\relax
    \fi
    \pgfkeysgetvalue{/pgfplots/cycle list shift}\pgfplots@listindex@shift
    \ifx\pgfplots@listindex@shift\pgfutil@empty
    \else
      \c@pgf@counta=\pgfplots@listindex\relax
      \advance\c@pgf@counta by\pgfplots@listindex@shift\relax
      \ifnum\c@pgf@counta<0
        \c@pgf@counta=-\c@pgf@counta
      \fi
      \edef\pgfplots@listindex{\the\c@pgf@counta}%
    \fi
    \ifnum\pgfplots@cycle@dim>0
      \c@pgf@counta=\pgfplots@cycle@dim\relax
      \c@pgf@countb=\pgfplots@listindex\relax
      \advance\c@pgf@counta by-1
      \pgfplotsloop{%
        \ifnum\c@pgf@counta<0
          \pgfplotsloopcontinuefalse
        \else
          \pgfplotsloopcontinuetrue
        \fi
      }{%
        \pgfkeysgetvalue{/pgfplots/cycle multi list/@N\the\c@pgf@counta}\pgfplots@cycle@N
        \pgfplotsmathmodint{\c@pgf@countb}{\pgfplots@cycle@N}%
        \divide\c@pgf@countb by \pgfplots@cycle@N\relax
        \expandafter\pgfplots@getautoplotspec@
        \csname pgfp@cyclist@/pgfplots/cycle multi list/@list\the\c@pgf@counta @\endcsname
        {\pgfplots@cycle@N}%
        {\pgfmathresult}%
        \t@pgfplots@toka=\expandafter{#1,}%
        \t@pgfplots@tokb=\expandafter{\pgfplotsretval}%
        \edef#1{\the\t@pgfplots@toka\the\t@pgfplots@tokb}%
        \advance\c@pgf@counta by-1
      }%
    \else
      \pgfplotslistsize\autoplotspeclist\to\c@pgf@countd

      \pgfplots@getautoplotspec@{\autoplotspeclist}{\c@pgf@countd}{\pgfplots@listindex}%
      \let#1=\pgfplotsretval
    \fi
    \pgfmath@smuggleone#1%
    \endgroup
  }

\pgfplotsset{
  cycle list set/.initial=
}
\makeatother

\usepackage{lipsum}
\usepackage{amsfonts}
\usepackage{graphicx}
\usepackage{epstopdf}
\usepackage{amsmath,amssymb}
\ifpdf
  \DeclareGraphicsExtensions{.eps,.pdf,.png,.jpg}
\else
  \DeclareGraphicsExtensions{.eps}
\fi

\headers{p-Robustness at Jacobi Speeds}{M. Wichrowski}

\title[Multigrid p-Robustness at Jacobi Speeds]{Multigrid p-Robustness at Jacobi Speeds: Efficient Matrix-Free Implementation of Local p-Multigrid Solvers}

\usepackage{amsopn}

\makeatletter
\newcommand*{\addFileDependency}[1]{
  \typeout{(#1)}
  \@addtofilelist{#1}
  \IfFileExists{#1}{}{\typeout{No file #1.}}
}
\makeatother

\setcopyright{none}
\renewcommand\footnotetextcopyrightpermission[1]{}

\author{Michał Wichrowski}
\affiliation{%
  \institution{Heidelberg University}
  \department{Interdisciplinary Center for Scientific Computing}
  \city{Heidelberg}
  \country{Germany}}
\email{mt.wichrowsk@uw.edu.pl}


\begin{document}

\begin{abstract}
  
Vertex-patch smoothers are essential for the robust convergence of geometric multigrid methods in high-order finite
element applications, yet their adoption is traditionally hindered by the prohibitive cost of solving local patch
problems. This paper presents a high-performance, matrix-free implementation of a p-multigrid local solver that
dismantles the trade-off between smoothing effectiveness and computational efficiency. We focus on the practical
realization of this iterative approach, leveraging sum-factorization and explicit SIMD vectorization to minimize memory
footprint and maximize arithmetic throughput. The performance analysis demonstrates that the solver effectively hides
data-fetching latencies and maintains optimal $\mathcal{O}(p^d)$ memory scaling, even when dominated by geometric data
on distorted meshes. The result is a robust smoother that rivals the execution speed of simple pointwise smoothers
while preserving the convergence benefits of patch-based methods.
\end{abstract}



\maketitle

\section{Introduction}

High-order finite element methods are a cornerstone of high-performance scientific computing, offering the potential
for rapid convergence and high accuracy. However, realizing this potential is contingent on the availability of
efficient and scalable solvers. Geometric multigrid methods with patch-based smoothers have emerged as a particularly
effective approach, demonstrating remarkable robustness with respect to the polynomial
degree~\cite{pavarino1993additive} and excellent data locality~\cite{wichrowski2025smoothers}, making them well-suited
for modern hardware.

The adoption of robust patch smoothers is hindered by two critical bottlenecks: the prohibitive computational cost of
the local solve and the substantial implementation effort required, particularly for complex physical models. This
creates a stark choice between the superior convergence of expensive vertex-patch methods and the raw speed of weak
cell-wise smoothers. We dismantle this trade-off. By engineering a fully matrix-free, p-multigrid local solver, we
achieve the mathematical robustness of Schwarz methods at the computational cost of a simple Jacobi iteration.

The development of efficient multigrid solvers is tightly linked to matrix-free computation and data locality. This
reflects a broader shift toward matrix-free techniques that are essential for high-order finite element
methods~\cite{Kronbichler2012,Kronbichler2017a, kronbichler2019multigrid} as traditional assembled sparse-matrix
approaches become prohibitively expensive due to massive memory requirements. Matrix-free methods circumvent this
bottleneck by never forming the global matrix, instead evaluating the operator's action on-the-fly. For methods on
tensor-product elements, sum-factorization techniques reduce the operator application cost to $\mathcal{O}(p^{d+1})$,
which is not only memory-optimal but also asymptotically cheaper than a conventional sparse matrix-vector
product~\cite{Kronbichler2017a}.

Since high-order methods are often memory-bound, performance hinges on efficient data movement. With modern hardware
characterized by deep memory hierarchies and varying access latencies, strategies to hide these latencies—such as
software prefetching~\cite{mahling2025fetch} or cache blocking—are essential. Patch-based approaches naturally excel in
this regard, as they decompose the global problem into smaller, local tasks that can fit into fast cache memory. This
concept aligns with segmental refinement strategies~\cite{adams2016segmental}, which demonstrate that processing local
subdomains can eliminate communication on fine grids and enhance locality. Recent work has shown that vertex-patch
smoothers provide excellent data locality specifically for high-order methods~\cite{wichrowski2025smoothers}, enabling
highly performant matrix-free implementations on both CPUs~\cite{munch2023cache, kronbichler2023enhancing} and
GPUs~\cite{cui2025implementation}. Our approach inherits these favorable properties, designed to operate entirely
within this matrix-free paradigm to preserve benefits at both the global and local levels. Furthermore, if coupled with
recent automatic differentiation techniques that streamline operator evaluation~\cite{wichrowski2025finitestrain}, the
framework presented here opens a path towards high-performance smoothers that are both computationally efficient and
easy to implement for complex problems.

While matrix-free techniques provide an efficient global framework, the strategy for solving the local patch problems
remains a critical design choice with significant trade-offs. \emph{Vertex-patch smoothers} function as overlapping
domain decomposition methods~\cite{pavarino1993additive}, offering remarkable robustness for high-order
discretizations~\cite{hong2016robust, KanschatMao15, voronin2025monolithic} and enabling advanced methods on
non-matching grids, such as CutFEM~\cite{bergbauer2025high, cui2025multigrid} or the Shifted Boundary
Method~\cite{wichrowski2025geometric, wichrowski2025matrix}. However, solving these patches often relies on
sophisticated direct methods, such as fast diagonalization~\cite{WitteArndtKanschat21} or static
condensation~\cite{brubeck2021scalable}, which can be inflexible or limited to structured grids. Conversely,
\emph{cell-wise smoothers}~\cite{margenberg2025hp, anselmann2024energy} are favored for their simplicity and low memory
footprint~\cite{fischer2000overlapping, pazner2018approximate, remacle2016gpu}, but they suffer from reduced smoothing
effectiveness.

Iterative methods align naturally with the matrix-free philosophy by avoiding assembly entirely, and such approaches
have been applied successfully to cell-based smoothers~\cite{bastian2019matrix} or discontinuous Galerkin
methods~\cite{fehn2020hybrid}. Following this rationale, a promising matrix-free and iterative solution was proposed
in~\cite{wichrowski2025local}, which employs a p-multigrid solver for the local patch problems. While that work
established the mathematical convergence of the method, this paper focuses on the software complexity, efficient
vectorization strategies, and performance optimization required to make this method competitive on modern hardware. Our
goal is to bridge the gap between algorithmic theory and high-performance execution, addressing the practical aspects
of the p-multigrid patch smoother, including its computational complexity, performance characteristics, and key
implementation details.

We validate the method on the Laplacian, but our primary contribution is a flexible and efficient framework for the
local patch-solve problem. The p-multigrid solver presented here is highly effective, yet modular enough to be adapted
or replaced by other advanced methods~\cite{brubeck2021scalable, WitteArndtKanschat21}, and the concept of a recursive
multigrid-based local solver is generalizable to more complex systems~\cite{braess1997efficient,wichrowski2022matrix,
  jodlbauer2024matrix}. We present and evaluate a practical, high-performance smoother implementation using the
open-source finite element library \texttt{deal.II}~\cite{dealII97,dealii2019design}. Key contributions include a
detailed description of a fully matrix-free implementation based on sum-factorization and the introduction of a
\emph{half V-cycle} optimization for the local solver, which significantly reduces computational cost by eliminating
the post-smoothing step. Finally, we provide a comprehensive performance evaluation, utilizing micro-benchmarking of
computational kernels to confirm optimal memory scaling of $\mathcal{O}(p^d)$ and to quantify the overhead from data
fetching operations inherent to patch-based methods.

Section~\ref{sec:multigrid} provides the mathematical background on multigrid methods, vertex patch smoothers, and the
specific p-multigrid local solver used in this work. Section~\ref{sec:implementation} details our matrix-free
implementation, focusing on data structures, vectorization strategies, and the optimization of the local solver.
Finally, Section~\ref{sec:performance} presents a comprehensive performance analysis, including micro-benchmarks of
computational kernels and scaling results on multi-core architectures.

\section{Mathematical and Algorithmic Background}
\label{sec:multigrid}
Let us first briefly describe the mathematical setup of the multigrid method, first
assume a hierarchy of meshes
\begin{gather}
  \mesh_0 \sqsubset \mesh_1 \sqsubset \dots \sqsubset \mesh_L,
\end{gather}
subdividing a domain in $\R^d$, where the symbol
``$\sqsubset$'' indicates nestedness, that is, every cell of mesh
$\mesh_{\ell+1}$ is obtained from a cell of mesh $\mesh_{\ell}$ by refinement. We
note that topological nestedness is sufficient from the algorithmic
point of view, such that domains with curved boundaries can be
covered approximately.

With each mesh $\mesh_\ell$ we associate a finite element space $V_\ell$ spanned by locally supported shape functions
with their associated degrees of freedom. We identify $V_\ell$ with $\mathbb R^{\operatorname{dim}V_\ell}$ and do not
distinguish between a finite element function $u_\ell$ and its coefficient vector by notation. Between these spaces, we
introduce transfer operators
\begin{gather}
  \begin{array}{rlcl}
    \resth\ell\colon & V_{\ell+1} & \to & V_\ell,     \\
    \prolh\ell\colon & V_{\ell}   & \to & V_{\ell+1}.
  \end{array}
\end{gather}
As usual, $\prolh\ell$ is chosen as the embedding operator and $\resth\ell$ as its
$\ell_2$-adjoint. The subscript $h$ indicates h-multigrid transfer between meshes of different resolution.

For our numerical experiments, we consider the Poisson equation with bilinear form
\begin{gather}
  \label{eq:bilinear-form}
  a(u,v) = \int_\Omega \nabla u \cdot \nabla v \, d\vec x.
\end{gather}
The discrete problem on level $\ell$ reads
\begin{gather}
  \label{eq:matrix}
  A_\ell u_\ell = b_\ell,
\end{gather}
where $A_\ell$ and $b_\ell$ are obtained by applying $a(\cdot,\cdot)$ and $b(\cdot)$ to the basis functions of
$V_\ell$.

Multigrid methods are common solution methods for the discrete linear system \eqref{eq:matrix}. The efficiency of the
overall method, however, depends on the choice of the so-called smoother, which is applied on each level in addition to
the operators $A_\ell$, $\resth\ell$, and $\prolh\ell$ which are rather standard. While simple smoothers like Jacobi
method are easy to implement, their performance often degrades when using high-order finite
elements~\cite{Kanschat08smoother, kronbichler2019multigrid, fehn2020hybrid, deville2002high}. In this work, we focus
on a vertex patch smoother, which has been shown to yield convergence rates that are independent of the polynomial
degree. We describe this smoother in detail in the next subsection.

\subsection{Vertex patch smoothers}
The vertex patch smoother is an overlapping subspace correction
method~\cite{ArnoldFalkWinther00,JanssenKanschat11,KanschatMao15,WitteArndtKanschat21,brubeck2021scalable,wichrowski2025smoothers}.
The domain is decomposed into a collection of overlapping subdomains, or \emph{patches} $\Omega_j$, each formed by the
cells surrounding an interior vertex. The smoother iteratively improves the solution by solving local problems on these
patches. We employ a multiplicative (Gauss-Seidel-like) variant, where patches are processed sequentially.

The update for a single patch $j$ consists of four steps. To define these, we introduce restriction operators:
$\overline{\Pi}_j$ extracts all degrees of freedom (DoFs) on a patch $\Omega_j$ (including its boundary) from a global
vector, while $\Pi_j$ extracts only the DoFs interior to the patch.
\begin{enumerate}
  \item \textbf{Gather.} Collect the current solution values on the patch into a local vector: $\overline{u}_j =
          \overline{\Pi}_j u$.
  \item \textbf{Evaluate.} Compute the local residual $r_j$ for the interior DoFs. This uses the local operator
        $\overline{A}_j$ (acting on all patch DoFs) and the global right-hand side $b$: $r_j = \Pi_j b - \Pi_j
          \overline{A}_j \overline{u}_j$.
  \item \textbf{Local solve.} Find a correction $d_j$ by approximately solving the local system on the patch
        interior: $d_j \approx \tilde A_j^{-1} r_j$, where $\tilde A_j^{-1}$ is an inexpensive approximate solver.
  \item \textbf{Scatter.} Add the local correction $d_j$ back to the global solution vector using the transpose of the
        interior restriction operator: $u \pluseq \Pi_j^T d_j$.
\end{enumerate}
These steps are summarized in Algorithm~\ref{alg:Loop-sequential}.  This implementation follows the cell-oriented loop
structure proposed in~\cite{wichrowski2025smoothers}, which is well-suited for matrix-free operator evaluation.

\begin{algorithm}[tp]
  \begin{algorithmic}
    \State \Function{\LUpNew}{j,u}
    \State    $ u_j \gets \overline\Pi_j u$ \Comment*{Gather}
    \State    $ r_{j} \gets \Pi_j b - \Pi_j\overline A_j u_j$ \Comment*{Evaluate}
    \State    $ d_j \gets \tilde{A}_j ^{-1} r_j$ \Comment*{Solve}
    \State \Return{ $\Pi_j^T d_j$ }\Comment*{Scatter}
    \EndFunction
    \State \For{$j=1,\ldots,N_\text{patches}$}	\Comment{Main smoother sweep}
    \State $u \pluseq$ \Call{\LUpNew}{$j$,$u$}
    \EndFor
  \end{algorithmic}
  \caption{Application of a patch smoother where local residuals are computed on-the-fly and combined with a local
    solver.}
  \label{alg:Loop-sequential}
\end{algorithm}

In this work, we focus on the local solve (step 3), which is critical for the overall performance. This step, denoted
as the application of an approximate inverse $\tilde{A}_j^{-1}$ in Algorithm~\ref{alg:Loop-sequential}, is performed by
applying a few iterations of a preconditioned Richardson method, as outlined in Algorithm~\ref{alg:local-richardson}.
Notice that the first iteration is performed outside the loop to avoid an unnecessary matrix-vector product, since the
initial guess for the correction is zero.

\begin{algorithm}[htbp]
  \caption{Local preconditioned Richardson iteration for the patch problem. This function implements the action of
  $\tilde{A}_j^{-1}$ from Algorithm~\ref{alg:Loop-sequential}.}
  \label{alg:local-richardson}
  \begin{algorithmic}[1]
    \Function{LocalSolve}{$\residual_j$}
    \Comment{The input $r_j$ is the residual for the patch problem.}
    \State $d_j \gets \Call{p-V-cycle}{\residual_j}$ \Comment{First iteration with zero initial guess.}
    \For{$k=2,\ldots,N_\text{iter}$} \Comment{Apply remaining steps.}
    \State $\res \gets \residual_j - A_j d_j$ \Comment{Compute residual for the correction.}
    \State $c_j \gets \Call{p-V-cycle}{\res}$ \Comment{Precondition with one p-MG cycle.}
    \State $d_j \pluseq c_j$ \Comment{Update correction.}
    \EndFor
    \State \Return{$d_j$}
    \EndFunction
  \end{algorithmic}
\end{algorithm}

\begin{figure}[tp]
  \centering
\def\svgwidth{\columnwidth}
\begin{tikzpicture}[scale=0.8, transform shape, font=\normalsize,
        cell/.style={draw=black, thick, minimum size=1.2cm},
        interior/.style={fill=blue!20},
        global/.style={fill=red!15},
        panel-label/.style={font=\bfseries\normalsize, align=center},
        arrow/.style={->, thick, shorten >=2mm, shorten <=2mm},
        dof/.style={circle, fill=black, inner sep=0.8pt}
    ]

    \newcommand{\drawPatchInteriorDoFs}[2]{
        \foreach \i in {1,...,5} {
                \foreach \j in {1,...,5} {
                        \pgfmathsetmacro{\xi}{\i/6}
                        \pgfmathsetmacro{\eta}{\j/6}

                        \coordinate (NW) at (#1.north west);
                        \coordinate (SE) at (#2.south east);

                        \coordinate (pt_x) at ($(NW)!\xi!(SE)$);	 
                        \coordinate (pt_y) at ($(NW)!\eta!(SE)$);  
                        \coordinate (pt) at (pt_x |- pt_y);
                        \node[dof, minimum size=4pt, inner sep=0pt, fill=blue] at (pt) {};
                    }
            }
    }

    \newcommand{\drawCubicDoFs}[1]{
        \foreach \xi in {0, 0.33, 0.66, 1} {
                \foreach \eta in {0, 0.33, 0.66, 1} {
                        \coordinate (bottom) at ($(#1.south west)!\xi!(#1.south east)$);
                        \coordinate (top) at ($(#1.north west)!\xi!(#1.north east)$);
                        \node[dof, minimum size=3pt, inner sep=0pt,  fill=red] at ($(bottom)!\eta!(top)$) {};
                    }
            }
    }

    \newcommand{\drawBoundaryDoFs}[2][dof]{
        \foreach \xi in {0, 0.33, 0.66, 1} {
                \coordinate (pt) at ($(#2.south west)!\xi!(#2.south east)$);
                \node[#1] at (pt) {};
            }
        \foreach \xi in {0, 0.33, 0.66, 1} {
                \coordinate (pt) at ($(#2.north west)!\xi!(#2.north east)$);
                \node[#1] at (pt) {};
            }
        \foreach \eta in {0.33, 0.66} {
                \coordinate (pt) at ($(#2.south west)!\eta!(#2.north west)$);
                \node[#1] at (pt) {};
            }
        \foreach \eta in {0.33, 0.66} {
                \coordinate (pt) at ($(#2.south east)!\eta!(#2.north east)$);
                \node[#1] at (pt) {};
            }
    }
    \newcommand{\drawInteriorFrame}[2]{%
        \draw[rounded corners=3pt, thick, draw=blue!60!black, fill=gray!80, fill opacity=0.1, densely dashed]
        ([shift={(4pt,-4pt)}]#1.north west)
        rectangle
        ([shift={(-4pt,4pt)}]#2.south east);%
    }

    \node[cell] (p1c1) at (0,0) {};
    \node[cell] (p1c2) at (1.2,0) {};
    \node[cell] (p1c3) at (0,-1.2) {};
    \node[cell] (p1c4) at (1.2,-1.2) {};
    \drawInteriorFrame{p1c1}{p1c4}

    \drawCubicDoFs{p1c1}
    \drawCubicDoFs{p1c2}
    \drawCubicDoFs{p1c3}
    \drawCubicDoFs{p1c4}
    \drawPatchInteriorDoFs{p1c1}{p1c4}

    \node[cell] (p2c1) at (4.5,0) {}; 
    \node[cell] (p2c2) at (5.7,0) {}; 
    \node[cell] (p2c3) at (4.5,-1.2) {}; 
    \node[cell] (p2c4) at (5.7,-1.2) {}; 
    \drawInteriorFrame{p2c1}{p2c4}
    \drawPatchInteriorDoFs{p2c1}{p2c4}

    \node[cell] (p3c1) at (9,0) {}; 
    \node[cell] (p3c2) at (10.2,0) {}; 
    \node[cell] (p3c3) at (9,-1.2) {}; 
    \node[cell] (p3c4) at (10.2,-1.2) {}; 
    \drawInteriorFrame{p3c1}{p3c4}
    \drawPatchInteriorDoFs{p3c1}{p3c4}

    \draw[arrow] ($(p1c3.west) - (2cm, 0)$) --	node[above, midway, ] {$\overline\Pi_j u, \;\;	\Pi_j b $} node[below,
        midway, ] {Gather}  (p1c3.west);
    \draw[arrow] (p1c4.east) -- node[above, midway, ] {$b_j - \Pi_j\overline A_j \overline u_j$} node[below, midway, ]
    {Evaluate} (p2c3.west);
    \draw[arrow] (p2c4.east) -- node[above, midway, ]{$\tilde A_j^{-1} r_j$}
    node[below, midway, ] {Local Solve}   (p3c3.west);
    \draw[arrow] (p3c4.east) -- node[above, midway, ]{$\Pi_j^T d_j$} node[below, midway, ] {Scatter} ($(p3c4.east) +
    (2cm, 0)$);

\end{tikzpicture}
  \caption{Workflow in a generic smoother application of a patch smoother for a single patch $j$. The three panels
    represent the main steps:
    gathering data from the global solution and right-hand side, computing the local residual, and applying the local
    correction and scattering it back to the global solution. Blue dots indicate interior DoFs of the patch, while red
    dots indicate all DoFs (including boundary) associated with each cell. The dashed rectangle highlights the patch
    interior.}
  \label{fig:patch-smoother-steps}
\end{figure}
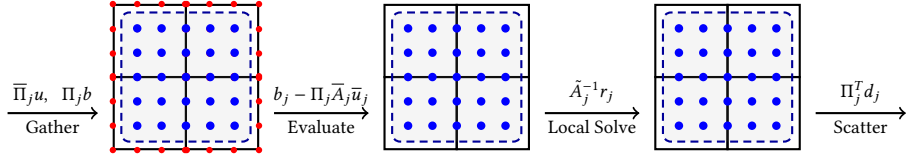

\subsection{The p-Multigrid Local Solver}
\label{sec:local-pMG-overview}
The efficiency of the patch smoother hinges on the preconditioner local solve step. To this end, we employ a multigrid
method. Since a patch contains too few cells for effective geometric coarsening, coarsening is instead performed by
reducing the polynomial degree $p$. This approach, detailed in~\cite{wichrowski2025local}, provides an effective and
computationally efficient way to approximate the action of the inverse. It is computationally inexpensive to apply
while keeping memory requirements low by relying on matrix-free data structures and storing only compact per-patch data
(diagonal entries).

As with any multigrid method, the p-multigrid solver is built upon a few key components: a hierarchy of levels,
transfer operators between these levels, a smoother and operator $A_l$ on each level, and a coarse-grid solver. In the
p-multigrid context, the levels correspond to a hierarchy of polynomial degrees. We choose a geometric progression,
such as $p_l=1, 3, 7, \dots$, up to the target degree. If the target degree does not fit this sequence, thehierarchy is
adjusted accordingly; for instance, for a target degree of $p=4$, the sequence would be $1, 3, 4$. The transfer
operators are the natural embedding for prolongation (from degree $p_{l-1}$ to $p_l$) and its adjoint for restriction.

The smoother on each p-level is a Richardson step with a preconditioner $P_{j,p}$. As proposed
in~\cite{wichrowski2025local}, we consider two choices for this preconditioner. The first is a damped Jacobi iteration,
where $P_{j,p}$ is the diagonal of the local operator $A_{j,p}$. This requires pre-computing and storing the inverse of
the diagonal for each patch and p-level. The second is a Cartesian-reinforced smoother, where the preconditioner
leverages the low cost of inverting the Laplacian on Cartesian grids using the fast diagonalization method. This
inverse is then combined with diagonal scaling to account for non-Cartesian geometry, yielding an efficient smoother
that outperforms Jacobi on nearly-Cartesian grids.

The complete p-multigrid V-cycle is outlined in Algorithm~\ref{alg:p-v-cycle} in which we also include the option for
skipping the post-smoothing step that we discuss in the next subsection. The algorithm recursively descends to the
coarsest level, where an exact solve is performed.

\begin{algorithm}[htbp]
  \caption{The p-multigrid V-cycle for a local patch problem. The post-smoothing step (lines 12-13) is omitted in the
    \emph{half V-cycle} optimization.}
  \label{alg:p-v-cycle}
  \begin{algorithmic}[1]
    \Function{V-cycle}{$l, \residual_p$}
    \If{$l = 1$}
    \State $d_1 \gets A_{j,1}^{-1} \residual_1$ \Comment{Coarse-level solve}
    \State \Return $d_1$
    \EndIf
    \State $d_p \gets \omega P_{j,p}^{-1} \residual_p$ \Comment{Pre-smoothing}
    \State $\res_p \gets \residual_p - A_{j,p} d_p$ \Comment{Compute residual}
    \State $\res_{p-1} \gets \restp{p,p-1} \res_p$ \Comment{Restrict residual}
    \State $e_{p-1} \gets \Call{V-cycle}{p-1, \res_{p-1}}$ \Comment{Recursive call}
    \State $e_p \gets \prolp{p-1,p} e_{p-1}$ \Comment{Prolongate correction}
    \State $d_p \gets d_p + e_p$ \Comment{Apply correction}
    \If{full V-cycle}
    \State $\res_p \gets \residual_p - A_{j,p} d_p$ \Comment{Compute residual for post-smoothing }
    \State $d_p \gets d_p + \omega P_{j,p}^{-1} \res_p$ \Comment{Post-smoothing }
    \EndIf
    \State \Return $d_p$
    \EndFunction
  \end{algorithmic}
\end{algorithm}

\subsubsection{The Half V-Cycle Optimization}
To reduce computational cost, we propose a \emph{half V-cycle} for the local p-multigrid solver. This is a standard
V-cycle but with the post-smoothing step omitted. This optimization saves one matrix-vector product and one diagonal
preconditioning step per V-cycle on each p-level, significantly reducing the cost of the local solve.

While this makes the local solver non-symmetric, the global multiplicative patch smoother is already non-symmetric due
to the fixed, sequential patch ordering. Therefore, we use GMRES as the outer solver, and this modification does not
change the choice of the global iterative method but reduces the computational work per smoother application. The half
V-cycle reduces the computational cost of the local solve at the expense of accuracy. A standard V-cycle requires three
local matrix-vector products (\texttt{vmult}s) per level, whereas the half V-cycle requires only two. However, since
the local solver starts with a zero initial guess for the correction, the \texttt{vmult} in the pre-smoothing step is
not required. This reduces the cost of the first full V-cycle to two \texttt{vmult}s at the highest polynomial degree,
and the cost of a half V-cycle to a single \texttt{vmult}. For any subsequent V-cycles, the cost would increase to
three and two \texttt{vmult}s, respectively. Since we apply only one local cycle, this optimization is significant.
Ultimately, the most important metric is the performance of the global solver. A cheaper, less accurate local solve is
beneficial if the total time to solution decreases. To investigate this trade-off, we compare the global GMRES
iteration counts when using a full V-cycle versus a half V-cycle as the local solver. We solve the Poisson problem on a
sequence of structured grids with increasing distortion to assess the robustness of each approach. To obtain distorted
mesh, each interior vertex of cartesian mesh $v$ is shifted by a distance $\delta h_v$ in a randomly chosen direction,
where $h_v$ is the minimum characteristic length of the edges attached to $v$ and $\delta$ is the distortion factor.
The results are presented in Table~\ref{tab:gmres_iterations_jacobi}.

\begin{table}[htbp]
  \centering
  \footnotesize
  \caption{GMRES iteration counts for a geometric multigrid preconditioner using a patch smoother with different local
    p-multigrid solver configurations. We compare a full V-cycle, a half V-cycle, and two half V-cycles against a
    nearly
    exact local solve (25 V-cycles).}
  \label{tab:gmres_iterations_jacobi}
  \begin{tabular}{l l |cccc|cccc}
    \toprule
      &                 & \multicolumn{4}{c|}{2D} & \multicolumn{4}{c}{3D}                     \\
    p & Local solver    & 0\%                     & 10\%                   & 25\% & 35\% & 0\%
      & 10\%            & 25\%                    & 30\%                                       \\
    \midrule
      & 1 V-cycle       & 5                       & 6                      & 7    & 9    & 6
      & 6               & 7                       & 7                                          \\
      & 1 half V-cycle  & 6                       & 7                      & 8    & 10   & 8
      & 9               & 10                      & 11                                         \\
    3 & 2 half V-cycles & 5                       & 6                      & 7    & 9    & 6
      & 6               & 7                       & 8                                          \\
    \myrowcolor
      & 25 V-cycles     & 4                       & 5                      & 6    & 7    & 4
      & 4               & 5                       & 5                                          \\
    \midrule
      & 1 V-cycle       & 5                       & 6                      & 8    & 10   & 8
      & 8               & 9                       & 9                                          \\
      & 1 half V-cycle  & 6                       & 7                      & 9    & 11   & 11
      & 11              & 13                      & 14                                         \\
    7 & 2 half V-cycles & 5                       & 7                      & 8    & 10   & 7
      & 8               & 9                       & 10                                         \\
    \myrowcolor
      & 25 V-cycles     & 3                       & 5                      & 6    & 7    & 3
      & 4               & 5                       & 5                                          \\
    \bottomrule
  \end{tabular}
\end{table}

Table~\ref{tab:gmres_iterations_jacobi} demonstrates the robustness of the proposed local solver configurations across
varying degrees of mesh distortion. While the nearly exact local solve (25 V-cycles) naturally yields the lowest
iteration counts, the approximate solvers remain highly effective. Specifically, the single half V-cycle incurs only a
modest penalty in the outer GMRES iterations compared to the full V-cycle, typically adding just one or two iterations
even on highly distorted grids. Interestingly, applying two half V-cycles often recovers the convergence behavior of
the full V-cycle or the exact solve, particularly for $p=3$. For the higher polynomial degree $p=7$, the single half
V-cycle shows slightly more sensitivity to strong distortions (e.g., 30\% in 3D).

\begin{figure}[tp]
  \input{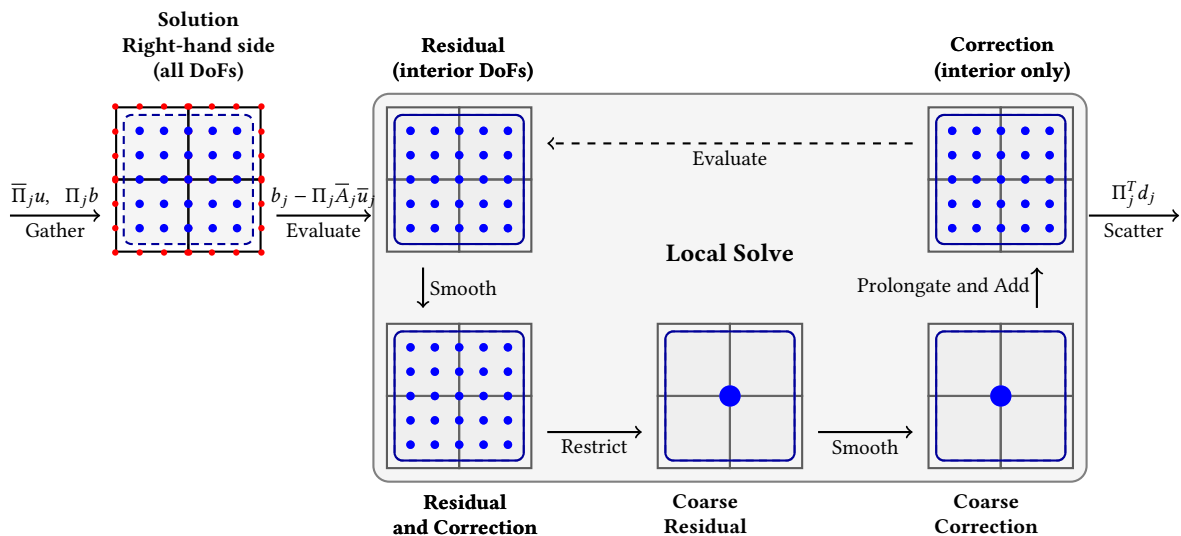}
  \caption{Patch-smoother data flow for a local p-multigrid solve (illustration for $p=3$ on the fine local level). The
    coarse level shown is $p=1$ and contains a single interior DoF (hence the coarse \emph{solve} is exact). The
    diagram
    shows gather, evaluate (local residual), pre-smoothing, restriction to the coarse level, an exact coarse
    correction, prolongation-and-add, and scatter-add. A dashed arrow indicates proceeding to the next local
    multigrid iteration. For brevity we illustrate only pre-smoothing; post-smoothing is omitted (\emph{half
      V-cycle}).}
  \label{fig:patch-flow}
\end{figure}

\section{The Implementation}
\label{sec:implementation}

The efficiency of high-order methods is tightly linked to matrix-free computation. Assembled sparse-matrix approaches
become prohibitively expensive as the polynomial degree $p$ increases, with memory requirements scaling as
$\mathcal{O}(p^{2d})$ and matrix-vector products as $\mathcal{O}(p^d)$. Matrix-free methods circumvent this by
evaluating the operator's action on-the-fly. For methods on tensor-product cells, sum-factorization techniques reduce
the operator application cost to $\mathcal{O}(p^{d+1})$, which is memory-optimal and asymptotically cheaper than a
sparse matrix-vector product~\cite{Kronbichler2017a}. Our implementation operates entirely within this paradigm, at
both the global and local levels.

Our implementation relies on the \texttt{deal.II} matrix-free framework, particularly the \texttt{FEEvaluation} class.
This framework is designed for high performance by processing cells in small batches, where the batch size is typically
chosen to match the hardware's SIMD vector width. Within each batch, cells are assigned a unique lane ID, enabling
vectorized computations across multiple cells simultaneously.

To use the \texttt{FEEvaluation} class, several pieces of data must be prepared for each cell batch. First, geometric
data, such as the Jacobian of the mapping from the reference cell, is computed. Second, the degrees of freedom
associated with each cell are identified. These tasks are automated by built-in functions within \texttt{deal.II}.

The framework internally detects and optimizes for common scenarios. For example, on meshes where cells are affine
images of the reference cell, the Jacobian is constant and can be efficiently cached. For Cartesian meshes, even faster
code paths are engaged, exploiting the fact that the Jacobian is diagonal or a scaled identity matrix. The framework
also identifies cells that share an identical mapping from the reference cell, computing and reusing geometric factors
for the entire group to minimize redundant calculations. These performance optimizations are handled transparently by
the \texttt{deal.II} matrix-free infrastructure, allowing the developer to focus on the physics of the problem.

As the core of evaluation we implement evaluation kernel as a separate class, that encapsulates the
quadrature-point-wise operations: computing contributions to the integral from provided gradients at quadrature points.
In case of our model problem \eqref{eq:bilinear-form} this is simply a pointwise multiplication of gradient with the
coefficient $\mu$. In more complex scenarios, such as finite-strain elasticity, this is where problem-defining material
parameters, like the elasticity tensor, would be stored~\cite{wichrowski2025finitestrain, davydov2021largeStrain}. The
functionality of \texttt{FEEvaluation} class then used to perform a vectorized evaluation of the operator's action on
each cell batch.

\subsection{Implementation of the Patch Smoother}
The practical implementation of Algorithm~\ref{alg:Loop-sequential} is realized using the \texttt{deal.II} finite
element library~\cite{dealII97,dealii2019design} and its matrix-free evaluation framework~\cite{Kronbichler2012,
  Kronbichler2017a}. This framework achieves high performance by processing cells in batches, utilizing SIMD
vectorization. Our implementation adapts this cell-based paradigm to the patch-based smoother. During an initialization
phase, we identify the cell indices that form each vertex patch.

Our implementation leverages SIMD vectorization by operating on the cells \emph{within} a single patch. In 2D, we
initialize a single \texttt{FEEvaluation} object with the list of cells that constitute the patch. This allows the
framework to assign each cell to a separate SIMD lane. For a quadrilateral mesh, a vertex patch consists of four cells.
This perfectly matches the vector width of AVX2 instruction sets, which operate on four double-precision numbers. Thus,
each of the four cells in the patch is processed in its own SIMD lane, allowing the operator evaluation to proceed in
parallel across the entire patch. In 3D, where a vertex patch on a hexahedral mesh consists of eight cells, we use two
\texttt{FEEvaluation} objects, providing a total of eight SIMD lanes to process all cells in two separate batches.

In addition to the storage inside \texttt{FEEvaluation} objects, we maintain local data buffers for storing the
gathered right-hand side, residual and correction vectors. These buffers are sized for the interior DoFs of a patch and
reused across patches to minimize memory allocation overhead.

The generic application of the smoother to a single patch $j$ involves the following detailed steps:
\begin{enumerate}
  \item \textbf{Re-initialize.} Re-initialize the cell-based evaluation kernels (\texttt{FEEvaluation}) with the list
        of cell IDs that constitute the current patch $j$. In this step, the data for the cells in
        the patch is fetched, which includes the degree of freedom indices and the Jacobians of the mapping from
        the reference cell and equation-specific data at quadrature points (e.g., the coefficient $\mu$ in
        \eqref{eq:bilinear-form}).
  \item \textbf{Fetch data.} Fetch the data required by operator evaluation and local solver. This includes the
        values
        of coefficient $\mu$ at quadrature points and the diagonal of the local operator $A_{j,p}$ for each
        p-level,
        which is needed for the Jacobi smoother as discussed in Section~\ref{sec:local-pMG-overview}.
  \item \textbf{Gather.} Collect the values of degrees of freedom for all cells within the patch from the global
        solution vector $u$ and the global right-hand side $b$ into cell-wise local buffers $\overline{u}_j$ and
        $\overline{b}_j$, respectively: $\overline{u}_j = \overline{\Pi}_j u$, $\overline{b}_j = \overline{\Pi}_j b$.
        This corresponds to applying the operator $\overline{\Pi}_j$ to both $u$ and $b$.
  \item \textbf{Cell-wise operator application.} For each cell in the patch, apply the local operator to the
        gathered solution values $\overline{u}_j$. This matrix-free evaluation uses the \texttt{FEEvaluation} objects
        for the patch, leveraging sum-factorization with a cost of $\mathcal{O}(p^{d+1})$ per cell, and produces a
        set of cell-local output vectors.
  \item \textbf{Local Gather.} The cell-local output vectors from the operator application are collected into a
        single patch-interior vector. The residual $r_j$ is then formed by subtracting this assembled vector from the
        gathered patch-interior right-hand side values.
  \item \textbf{Local solve.} Apply the local solver to the local residual $r_j$ to obtain the correction $d_j$. This
        step is performed entirely on the \emph{local} data buffers.
  \item \textbf{Scatter-add.} Add the computed correction $d_j$ from the local buffer back into the global solution
        vector $u$, corresponding to the action of $\Pi_j^T$.
\end{enumerate}
The final scatter-add step is implemented in two stages to reuse the existing cell-based infrastructure of the
\texttt{FEEvaluation} class. First, the correction vector $d_j$, which is organized by patch-interior degrees of
freedom, is distributed into a temporary cell-wise local buffer.
To correctly handle degrees of freedom shared by multiple cells within the patch, this operation
is masked. For each shared DoF, the correction is transferred only to the first cell that contains it; contributions
from subsequent cells are masked out to prevent duplicate updates. Second, the \texttt{FEEvaluation} object's built-in
scatter functionality is invoked to add the values from this local buffer to the corresponding entries in the global
solution vector $u$.

\subsection{Implementation of the Local  p-Multigrid Solver}
We focus on the implementation of the local solver, which we view as an iterative method, even though it may be applied
for only a single iteration. As the solver we use a local p-multigrid method, which requires the standard multigrid
components: a matrix-vector product, inter-grid transfer operators, and a coarse-grid solver. On top of these, a
smoother is needed on each p-level. For efficiency, the \texttt{deal.II} \texttt{FEEvaluation} framework is templated
on the polynomial degree, which allows the compiler to unroll loops and generate highly optimized
code~\cite{Kronbichler2017a}. We preserve this design principle in our local solver; the p-multigrid functions are
templated on the local p-level, which in turn defines the polynomial degree for the operators on that level.
\begin{figure}[tp]
  \centering
  \input{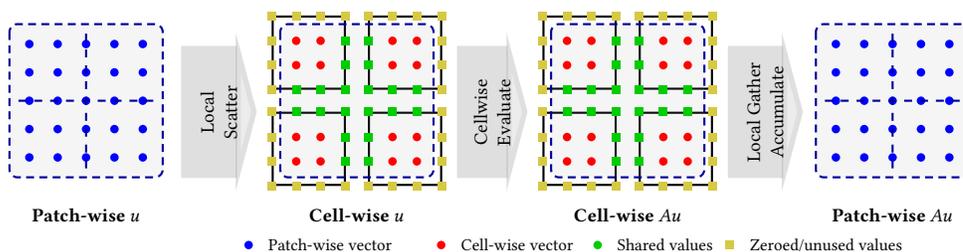}
  \caption{Application of vmult in the patch-local solver }
  \label{fig:vmult-patch}
\end{figure}

The smoother within the local p-multigrid cycle is a damped Jacobi iteration. Its implementation is straightforward: it
involves a multiplication of the residual vector by the inverse of the operator's diagonal, scaled by a damping factor
$\omega$. The inverses of the diagonal for the local operator $A_{j,p}$ must be pre-computed and stored for each patch
$j$ and for each polynomial degree $p$ used in the local p-multigrid hierarchy.

Alternatively, in~\cite{wichrowski2025local} a Cartesian-reinforced smoother was used which exploits the low-cost of
application of inverse of Laplacian on Cartesian grids via fast diagonalization technique. By combining the inverse
obtained via fast diagonalization with diagonal scaling to account for non-Cartesian geometry, an efficient smoother is
obtained. As expected for nearly Cartesian grids, this smoother outperforms the Jacobi smoother.

The transfer operators, which prolongate corrections from a coarser p-level ($p-1$) to a finer one ($p$) and restrict
residuals in the opposite direction, are also implemented efficiently. For tensor-product elements, these operators
have a Kronecker product structure composed of 1D transfer operators. We use a pre-computed 1D transfer matrix and
apply it sequentially along each coordinate direction. For a 3D problem, this corresponds to applying operators of the
form $I \otimes I \otimes T$, then $I \otimes T \otimes I$, and finally $T \otimes I \otimes I$, where $T$ is the 1D
transfer matrix and $I$ is the identity. This results in a computational complexity of $\mathcal{O}(p^{d+1})$. A
significant part of 1D transfer matrices is composed of zeros which we exploit to reduce the number of arithmetic
operations. The implementation consists of simple loops whose bounds are known at compile time due to the templating on
the p-level allowing the compiler to unroll and optimize them.

The coarse-grid problem in the p-multigrid hierarchy corresponds to the local patch problem with linear finite elements
($p=1$). For quadrilateral or hexahedral meshes, a patch has only a single degree of freedom in its interior at this
level --- the one associated with the central vertex itself. Consequently, the local interior system $A_{j,1}$ is a $1
  \times 1$ matrix. The coarse-grid \emph{solve} thus simplifies to a single scalar multiplication, making it both exact
and computationally trivial.

\subsubsection{Local vmult}
The core computational kernel of the local solver is the matrix-vector product, which we refer to by its common
function name \texttt{vmult}. This operation computes the action of the local operator $A_{j,p}$ on a vector defined on
the patch interior. The \texttt{vmult} implementation follows the cell-based paradigm of the matrix-free framework. It
first takes the input vector, defined on the patch's interior DoFs, and distributes its values to the cell-local data
structures used by the \texttt{FEEvaluation} objects. Then, each \texttt{FEEvaluation} object performs the operator
evaluation on its assigned cells using sum-factorization. Finally, the contributions from all cells in the patch are
gathered and summed into the output patch vector. This process is illustrated in Figure~\ref{fig:vmult-patch}. This
design inherits the optimal computational complexity of $\mathcal{O}(p^{d+1})$ and leverages vectorization across the
cells of the patch.

A critical aspect of the local \texttt{vmult} is the data reshuffling required to interface between the patch-level
data layout and the cell-based evaluation framework. The input vector, organized by patch-interior DoFs, must first be
distributed to the cell-local data structures expected by the \texttt{FEEvaluation} objects. After the cell-wise
operator evaluation, the resulting local vectors must be gathered and summed back into a single patch-level output
vector. This data reshuffling occurs twice for every matrix-vector product. We consider two strategies for managing the
index mapping required for these distribute and gather operations: either pre-computing the mappings into compile-time
look-up tables, which offers fast access at the cost of memory, or dynamically computing the indices on-the-fly, which
reduces the instruction count but may introduce computational overhead.

To choose the more performant indexing strategy, we compare the performance of the two approaches discussed:
pre-computing the index mappings versus computing them on-the-fly. The on-the-fly computation of indices consistently
outperformed the pre-computed look-up table approach in our benchmarks. This is likely due to the simple, arithmetic
nature of the index calculation, which can be heavily optimized by the compiler and benefits from better instruction
cache locality compared to the memory accesses required for a look-up table. Interestingly, we observed that the
\texttt{clang} compiler (version 15) produced significantly more efficient code for this particular kernel than
\texttt{gcc} (version 11.2), leading to the best overall performance.

A notable aspect of this method is that the entire smoother application, including the recursive p-multigrid V-cycle,
can be executed without a single floating-point division. During the initialization phase, we compute and store the
inverse of the diagonal for each p-level, $D_{j,p}^{-1}$. Consequently, the coarse-grid solve and the Jacobi smoothing
steps are transformed into multiplications. This avoids more expensive division instructions during the main solver
loop, allowing the core computations to be expressed entirely through additions and multiplications, which can be
efficiently mapped to fused multiply-add (FMA) instructions on modern hardware.

The primary memory cost associated with the smoother comes from storing the diagonal of the local operator $A_{j,p}$
for each patch $j$ and p-level $p$, which is needed for the Jacobi smoother. Since we only store the diagonal, the
memory per degree of freedom is constant, and the total storage for these diagonals scales optimally as
$\mathcal{O}(p^d)$ per patch.

\section{Performance Analysis}
\label{sec:performance}
To assess the performance under different geometric conditions, we consider two representative test cases. The first is
a standard Cartesian grid, where all cells are identical axis-aligned boxes. The second is a distorted grid, generated
by randomly perturbing the vertices of a Cartesian grid by up to 10\% of the cell size. The primary difference between
these scenarios lies in the structure of the Jacobian matrix for the mapping from the reference cell to a physical
cell. For the numerical experiments we build the mesh hierarchy by starting from a single patch (one coarse cell) and
applying uniform global refinements to obtain the finer levels used in the benchmarks. To obtain distorted mesh, each interior vertex $v$ is shifted by a distance $\delta h_v$ in a randomly chosen direction, where
$h_v$ is the minimum characteristic length of the edges attached to $v$ and $\delta$ is the distortion factor.

For our numerical experiments, we select a range of polynomial degrees: $p=2, 3, 4,$ and $7$. The lower degrees are
chosen as they represent common and practical choices in many finite element applications. The choice of $p=7$ is
motivated by its alignment with the geometric coarsening sequence of our p-multigrid solver, which naturally steps down
from $p=7$ to $p=3$ and then to $p=1$. Furthermore, $p=7$ corresponds to 8 degrees of freedom per spatial direction,
resulting in a total number of cell DoFs that is a power of two ($8^d$). While our current implementation does not
explicitly optimize for this property, it presents an interesting avenue for future hardware-specific optimizations.

For the Cartesian grid, the Jacobian is diagonal, which allows for highly optimized evaluation kernels with minimal
data movement. In contrast, for the distorted grid, the Jacobian is a full matrix that varies between quadrature
points. This necessitates storing and fetching all $d \times d$ components of the inverse Jacobian for each quadrature
point, significantly increasing the memory bandwidth pressure and the number of floating-point operations during
operator evaluation. Additionally, this complexity impacts the data reshuffling phase when fetching patch-wise
geometrical data, as the process must loop over all quadrature points to gather the varying geometric factors.
Comparing these two cases allows us to quantify the performance impact of geometric complexity.

We begin our performance analysis with a series of micro-benchmarks designed to dissect the computational cost of the
patch smoother. To isolate the algorithmic efficiency of the local solver components from main memory latency, we
employ a small test problem consisting of 343 patches on the finest level. This problem size is sufficiently small
(Table~\ref{tab:dofs-small}) to fit entirely within the CPU cache, ensuring that our measurements, performed on a
single core, reflect the computational behavior of the algorithms. The following results were obtained on a dual-socket
node equipped with AMD EPYC 7282 processors, providing a total of 32 physical cores.

\begin{table}[htbp]
  \centering
  \caption{Problem sizes for microbenchmarking for different polynomial degree $p$.}
  \label{tab:dofs-small}
  \begin{tabular}{lrrrr}
    \toprule
    Patches & \multicolumn{4}{c}{Number of DoFs}
    \\                                                    & $p=2$ & $p=3$ & $p=4$ & $p=7$
    \\ \midrule
    342     & 4,913                              & 15,625 & 35,937 & 185,193 \\
    \bottomrule
  \end{tabular}
\end{table}

\subsection{Microbenchmarking of Local Solver Components}
\label{sec:benchmarking}

\begin{figure}[htbp]
  { \input{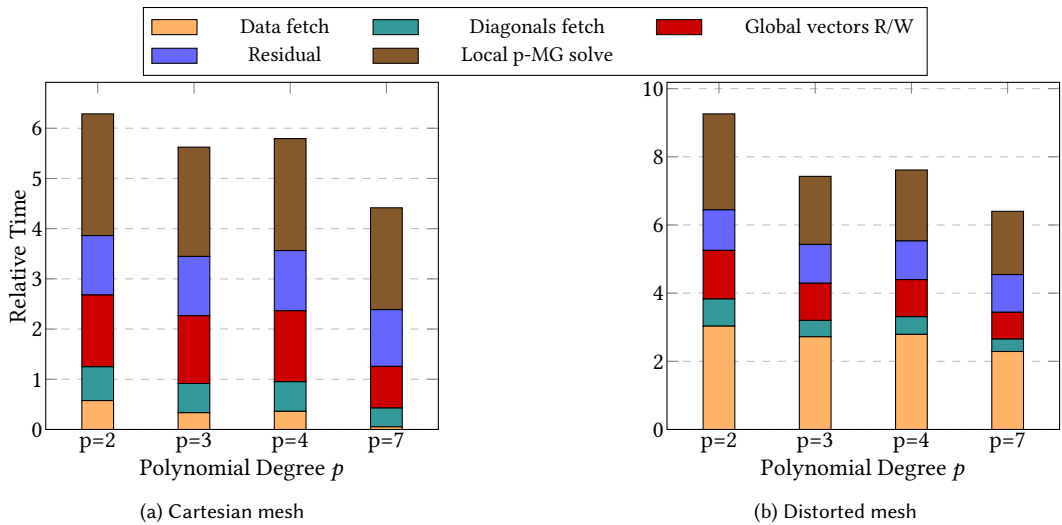} }
  \caption{Breakdown of the wall-clock time for a single patch smoother application across different polynomial
    degrees
    $p$. Left: Cartesian mesh. Right: distorted mesh. Components: initial data fetching and setup, local
    residual computation, the p-multigrid local solve, and scattering the correction back to the global vector.
    Relative times are normalized by the time required to evaluate the patch-local operator $A$ application.}
  \label{fig:patch-application-breakdown-subfigs}
\end{figure}

Figure~\ref{fig:patch-application-breakdown-subfigs} presents the breakdown of the wall-clock time for a single
smoother application, distinguishing between the Cartesian and distorted mesh cases. The reported times cover the key
phases: data fetching and setup, residual computation, the local p-multigrid solve, and the scattering of results. To
obtain reliable measurements for these short kernels, we record LIKWID performance counters~\cite{likwid} during the
actual smoother application. Each measurement consists of 20 timed repetitions of the multiplicative patch smoother.
For each repetition we record LIKWID performance counters and the wall-clock time. Reported timings are the arithmetic
mean over the 20 repetitions. All times are given relative to (i.e., normalized by) the time required to evaluate the
patch-local operator application $A u$ (a single \texttt{vmult}); thus a reported value of 2.0 indicates an operation
twice as expensive as one patch-wise $A$ evaluation. This normalization allows for a more consistent comparison across
different polynomial degrees, as the cost of the baseline operator scales with $p$. Furthermore, while performance
measurement tools like LIKWID may introduce a small overhead, reporting relative timings helps to mitigate this effect,
ensuring a fair comparison of the intrinsic cost of each component.

For the Cartesian mesh, the geometric mappings are constant, allowing for highly optimized data access. As shown in the
results (corresponding to the left panel of Figure~\ref{fig:patch-application-breakdown-subfigs}), the data fetching
phase is negligible when measured relative to the patch-local operator evaluation. The runtime is therefore dominated
by the arithmetic intensity of the solver components. For instance, at $p=2$, the local p-multigrid solve is
approximately 2.42 times the cost of a single patch-wise $A$ evaluation, while data fetching is about 0.57 times that
cost. Similarly, at $p=7$, the solve is about 2.03 times the single-$A$ cost and fetching is only about 0.052 times.
This confirms that on structured grids the method is compute-bound and effectively utilizes the matrix-free operator
evaluation.

The component labeled 'Global vectors R/W' captures the cost of accessing the global solution and right-hand side
vectors during the gather step, as well as updating the solution during the scatter step. Even on this small benchmark,
designed to minimize main memory latency, these operations constitute a noticeable fraction of the total runtime. This
is particularly evident on the Cartesian grid, where the geometric setup is minimal. Here, the cost of reading and
writing global vector entries is more than twice that of fetching the geometric data, highlighting that even efficient
indirect addressing into global vectors incurs a fixed overhead that cannot be fully optimized away.

The performance profile changes significantly for the distorted mesh (right panel of
Figure~\ref{fig:patch-application-breakdown-subfigs}). In this scenario, the mapping is general, necessitating the
storage and retrieval of full Jacobian matrices at each quadrature point. This geometric complexity imposes a heavy
penalty on the setup phase. The cost of fetching the data increases dramatically. At $p=7$, the data fetching time
rises to 2.28 times cost of a single patch-wise $A$ evaluation that is exceeding the time required for the local solve
itself (1.86 times $A$ evaluation). Even at lower polynomial degrees, such as $p=2$, the fetching cost (3.03 times $A$
evaluation) is comparable to the solve time (2.81 times $A u$ evaluation). This demonstrates that for complex
geometries, the overhead of managing varying geometric factors becomes a dominant factor, shifting the bottleneck from
arithmetic computation to data movement and setup.

\paragraph{CPU Throughput}
To understand the performance characteristics of the patch smoother, we analyze the computational throughput of its
constituent parts, measured in millions of floating-point operations per second (GFLOP/s). High throughput indicates
that the hardware's floating-point units are being used effectively, which is typical for arithmetically intense
computations. Conversely, low throughput often signals that the operation is limited by memory bandwidth, with the
processor spending most of its time waiting for data.

\begin{figure}[htbp]
  \centering

\newcommand{\ProcesspatchPatchBreakDown}[2]{%
    \pgfkeys{/pgf/fpu=true, /pgf/fpu/output format=fixed}%
    \pgfplotstabletranspose[
        colnames from=section,
        input colnames to=property
    ]\datatransposedtemp{#1}

    \pgfplotstablegetelem{0}{PATCH_compute_A_residual}\of{\datatransposedtemp}
    \let\computeresidual=\pgfplotsretval
    \pgfplotstablegetelem{0}{PATCH_Au_to_residual}\of{\datatransposedtemp} \let\autoresidual=\pgfplotsretval
    \pgfplotstablegetelem{0}{PATCH_local_solve}\of{\datatransposedtemp} \let\localsolve=\pgfplotsretval

    \expandafter\pgfmathsetmacro\csname #2AuThrough\endcsname{(\computeresidual/1000)}
    \expandafter\pgfmathsetmacro\csname #2SolveThrough\endcsname{\localsolve/1000}
    \expandafter\pgfmathsetmacro\csname #2ToResidualThrough\endcsname{\autoresidual/1000}
    \pgfkeys{/pgf/fpu=false}%
}

\newcommand{\ReadAndProcessPatchFile}[2]{%
    \pgfplotstableread[
        col sep=comma,
        header=true,
        comment chars={-},
        trim cells=true
    ]{#1}\loadedpatchdatatable
    \ProcesspatchPatchBreakDown{\loadedpatchdatatable}{#2}%
}


\begin{subfigure}[b]{0.49\textwidth}
    \centering

    \ReadAndProcessPatchFile{results/cartesian/small_p2.csv}{Two}

    \ReadAndProcessPatchFile{results/cartesian/small_p3.csv}{Three}

    \ReadAndProcessPatchFile{results/cartesian/small_p4.csv}{Four}
    \ReadAndProcessPatchFile{results/cartesian/small_p7.csv}{Seven}
    \begin{tikzpicture}

        \pgfplotstableread[col sep=comma]{
        label,p2,p3,p4,p7
        {Evaluate\\operator},\TwoAuThrough,\ThreeAuThrough,\FourAuThrough,\SevenAuThrough
        {p-MG\\solve},\TwoSolveThrough,\ThreeSolveThrough,\FourSolveThrough,\SevenSolveThrough

        {Local\\Gather},\TwoToResidualThrough,\ThreeToResidualThrough,\FourToResidualThrough,\SevenToResidualThrough
        }\datatablecartesian

        \begin{axis}[
                width=\linewidth,
                height=6.2cm,
                ybar,
                bar width=5pt,
                enlarge x limits=0.25,
                ylabel={Throughput [GFLOP/s]},
                symbolic x coords={{Evaluate\\operator}, {p-MG\\solve}, {Local\\Gather}},
                xtick=data,
                xticklabel style={align=center, font=\small},
                ymin=0, ymax=30,
                ymajorgrids=true,
                grid style=dashed,
                legend style={at={(0.5,1.05)}, anchor=south, legend columns=-1, font=\small},
            ]
            \addplot[fill=BarColorOne] table[x=label, y=p2] {\datatablecartesian};
            \addlegendentry{$p=2$}
            \addplot[fill=BarColorTwo] table[x=label, y=p3] {\datatablecartesian};
            \addlegendentry{$p=3$}
            \addplot[fill=BarColorThree] table[x=label, y=p4] {\datatablecartesian};
            \addlegendentry{$p=4$}
            \addplot[fill=BarColorFour] table[x=label, y=p7] {\datatablecartesian};
            \addlegendentry{$p=7$}
        \end{axis}
    \end{tikzpicture}
    \caption{Cartesian mesh}
    \label{fig:throughput-cartesian}
\end{subfigure}
\hfill
\begin{subfigure}[b]{0.49\textwidth}

    \ReadAndProcessPatchFile{results/cartesian/small_p2.csv}{Two}

    \ReadAndProcessPatchFile{results/cartesian/small_p3.csv}{Three}

    \ReadAndProcessPatchFile{results/cartesian/small_p4.csv}{Four}
    \ReadAndProcessPatchFile{results/cartesian/small_p7.csv}{Seven}
    \centering
    \begin{tikzpicture}
        \pgfplotstableread[col sep=comma]{
        label,p2,p3,p4,p7
        {Evaluate\\operator},\TwoAuThrough,\ThreeAuThrough,\FourAuThrough,\SevenAuThrough
        {p-MG\\solve},\TwoSolveThrough,\ThreeSolveThrough,\FourSolveThrough,\SevenSolveThrough

        {Local\\Gather},\TwoToResidualThrough,\ThreeToResidualThrough,\FourToResidualThrough,\SevenToResidualThrough
        }\datatabledistorted

        \begin{axis}[
                width=\linewidth,
                height=6.2cm,
                ybar,
                bar width=5pt,
                enlarge x limits=0.25,
                symbolic x coords={{Evaluate\\operator}, {p-MG\\solve}, {Local\\Gather}},
                xtick=data,
                xticklabel style={align=center, font=\small},
                ymin=0, ymax=30,
                ymajorgrids=true,
                grid style=dashed,
                legend style={at={(0.5,1.05)}, anchor=south, legend columns=-1, font=\small},
            ]
            \addplot[fill=BarColorOne] table[x=label, y=p2] {\datatabledistorted};
            \addlegendentry{$p=2$}
            \addplot[fill=BarColorTwo] table[x=label, y=p3] {\datatabledistorted};
            \addlegendentry{$p=3$}
            \addplot[fill=BarColorThree] table[x=label, y=p4] {\datatabledistorted};
            \addlegendentry{$p=4$}
            \addplot[fill=BarColorFour] table[x=label, y=p7] {\datatabledistorted};
            \addlegendentry{$p=7$}
        \end{axis}
    \end{tikzpicture}
    \caption{Distorted mesh}
    \label{fig:throughput-distorted}
\end{subfigure}
  \caption{Throughput of individual components of the patch smoother application. The y-axis shows
    performance in GFLOP/s.}
  \label{fig:throughput-breakdown}
\end{figure}
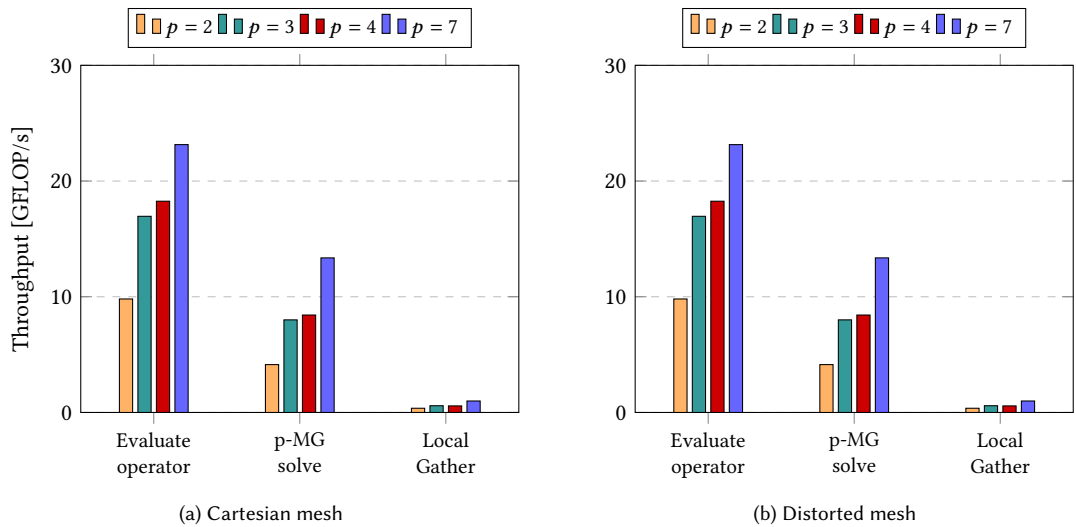

Figure~\ref{fig:throughput-breakdown} presents a breakdown of the throughput for each major step in a single patch
application. The results clearly distinguish between memory-bound and compute-bound operations.

The core computational kernels---the operator evaluation (`Evaluate operator`) and the local p-multigrid solve (`p-MG
solve`)---achieve high throughput, approaching the performance of a standard, highly optimized matrix-vector product.
The `Local Gather` operation, which assembles the patch residual from cell-local contributions, exhibits significantly
lower throughput, as it involves non-contiguous memory accesses and data reshuffling that limit its efficiency.

We omit the data fetching and scattering steps from the analysis as as they perform very few arithmetic operations
relative to the amount of data they move and thus their performance is order of magnitude lower. Operations like
fetching cell data and reading from global solution vectors are entirely memory-bound. The scatter operation
(`Distribute correction`) performs slightly better as it involves some arithmetic for accumulating values, but it is
still fundamentally limited by memory access patterns. This analysis highlights a critical performance challenge in
patch-based methods: while the local computations are efficient, the overhead of gathering and scattering data across
the non-contiguous memory layout of a patch can dominate the overall runtime.

\subsection{Breakdown of p-MG performance}
\begin{figure}[htbp]
  \input{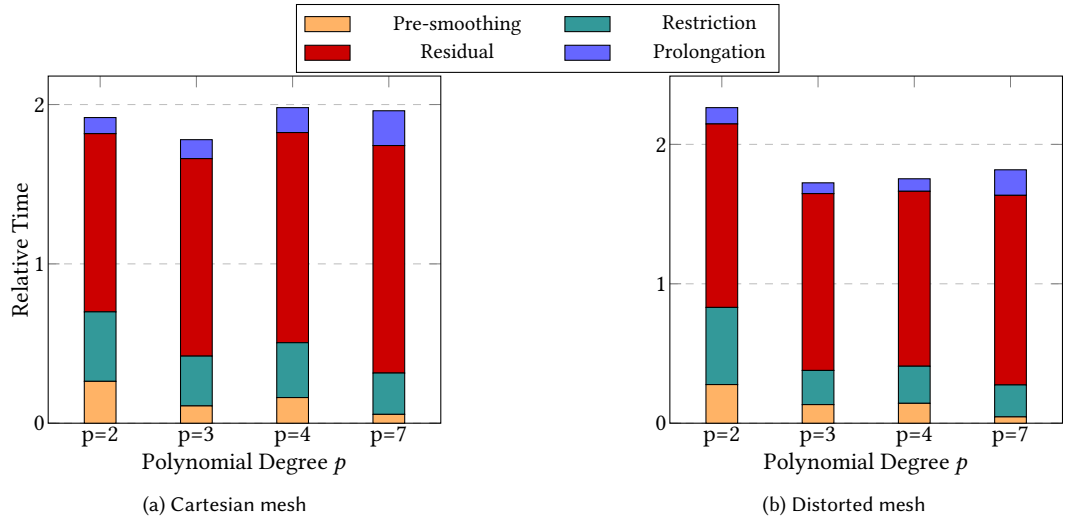}
  \caption{Breakdown of the wall-clock time for a single application of local p-multigrid across different polynomial
    degrees
    $p$. Left: Cartesian mesh. Right: distorted mesh. Components: residual evaluation, restriction,
    prolongation, and smoothing.}
  \label{fig:breakdown-pmg}
\end{figure}

Having analyzed the overhead of the outer patch loop, we now turn our attention to the internal performance of the
local p-multigrid solver itself. Figure~\ref{fig:breakdown-pmg} provides a detailed breakdown of the wall-clock time
spent in the four primary components of the p-multigrid V-cycle: pre-smoothing, residual computation, restriction, and
prolongation. As before, these timings are normalized by the cost of a single fine-level patch-wise operator evaluation
($Au$). It is important to note that the reported times are cumulative, representing the sum of time spent in each
component across all levels of the p-multigrid hierarchy.

The results indicate that the residual computation is the dominant cost, consistently exceeding the time of a single
fine-level operator evaluation (a relative time greater than 1.0). This is expected for several reasons. First, the
residual calculation involves a matrix-vector product ($Au$) on each level of the hierarchy; while the cost decreases
rapidly with $p$, the aggregate cost is naturally higher than the fine-level evaluation alone. Second, the residual
step includes additional vector arithmetic for subtraction and, crucially, the overhead of reshuffling data between the
formats required for the operator evaluation and the level-specific vectors.

The inter-level transfer operators (restriction and prolongation) consume approximately one-quarter of the time
required for the local operator evaluation. This performance aligns with theoretical expectations, as these operators
share the $\mathcal{O}(p^{d+1})$ complexity of the matrix-vector product but involve simpler arithmetic structures. The
smoothing step, implemented as a damped Jacobi iteration, is the least expensive component. Its cost remains relatively
constant across polynomial degrees, reflecting its linear complexity with respect to the number of degrees of freedom
on the patch.

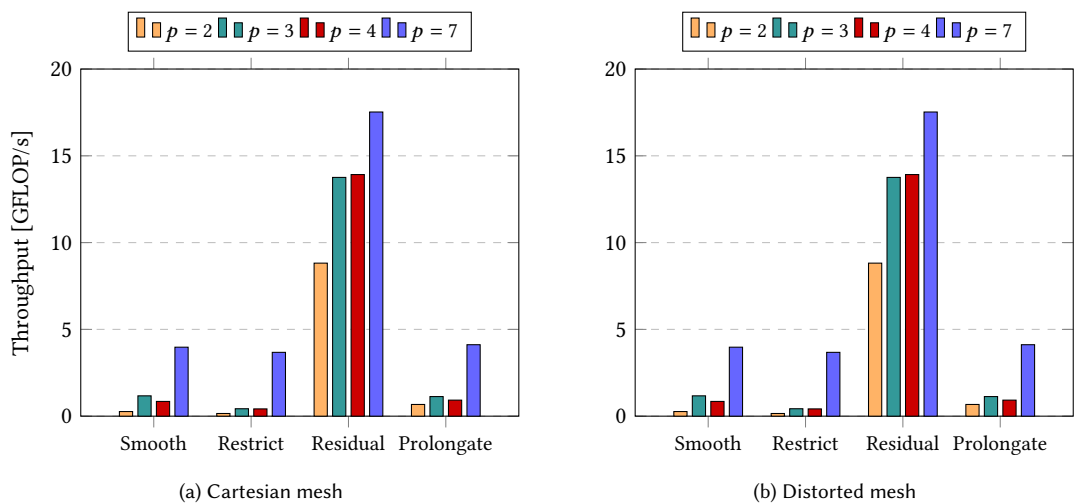
\begin{figure}[htbp]
  \centering

\newcommand{\ProcesspatchPatchBreakDown}[2]{%
    \pgfkeys{/pgf/fpu=true, /pgf/fpu/output format=fixed}%
    \pgfplotstabletranspose[
        colnames from=section,
        input colnames to=property
    ]\datatransposedtemp{#1}

    \pgfplotstablegetelem{0}{MG_Pre-smooth}\of{\datatransposedtemp}
    \let\presmooth=\pgfplotsretval
    \pgfplotstablegetelem{0}{MG_Residual}\of{\datatransposedtemp} \let\residual=\pgfplotsretval
    \pgfplotstablegetelem{0}{MG_Restrict}\of{\datatransposedtemp} \let\restrict=\pgfplotsretval
    \pgfplotstablegetelem{0}{MG_Prolongate}\of{\datatransposedtemp} \let\prologate=\pgfplotsretval

    \expandafter\pgfmathsetmacro\csname #2Presmooth\endcsname{(\presmooth/1000)}
    \expandafter\pgfmathsetmacro\csname #2Restrict\endcsname{\restrict/1000}
    \expandafter\pgfmathsetmacro\csname #2Residual\endcsname{\residual/1000}
    \expandafter\pgfmathsetmacro\csname #2ProlongateThrough\endcsname{\prologate/1000}
    \pgfkeys{/pgf/fpu=false}%
}

\newcommand{\ReadAndProcessPatchFile}[2]{%
    \pgfplotstableread[
        col sep=comma,
        header=true,
        comment chars={-},
        trim cells=true
    ]{#1}\loadedpatchdatatable
    \ProcesspatchPatchBreakDown{\loadedpatchdatatable}{#2}%
}


\begin{subfigure}[b]{0.49\textwidth}
    \centering

    \ReadAndProcessPatchFile{results/cartesian/pmg_p2.csv}{Two}

    \ReadAndProcessPatchFile{results/cartesian/pmg_p3.csv}{Three}

    \ReadAndProcessPatchFile{results/cartesian/pmg_p4.csv}{Four}
    \ReadAndProcessPatchFile{results/cartesian/pmg_p7.csv}{Seven}
    \begin{tikzpicture}

        \pgfplotstableread[col sep=comma]{
        label,p2,p3,p4,p7
        {Smooth},\TwoPresmooth,\ThreePresmooth,\FourPresmooth,\SevenPresmooth
        {Restrict},\TwoRestrict,\ThreeRestrict,\FourRestrict,\SevenRestrict

        {Residual},\TwoResidual,\ThreeResidual,\FourResidual,\SevenResidual
        {Prolongate},\TwoProlongateThrough,\ThreeProlongateThrough,\FourProlongateThrough,\SevenProlongateThrough
        }\datatablecartesian

        \begin{axis}[
                width=\linewidth,
                height=6.2cm,
                ybar,
                bar width=5pt,
                enlarge x limits=0.25,
                ylabel={Throughput [GFLOP/s]},
                symbolic x coords={{Smooth}, {Restrict}, {Residual}, {Prolongate}},
                xtick=data,
                xticklabel style={align=center, font=\small},
                ymin=0, ymax=20,
                ymajorgrids=true,
                grid style=dashed,
                legend style={at={(0.5,1.05)}, anchor=south, legend columns=-1, font=\small},
            ]
            \addplot[fill=BarColorOne] table[x=label, y=p2] {\datatablecartesian};
            \addlegendentry{$p=2$}
            \addplot[fill=BarColorTwo] table[x=label, y=p3] {\datatablecartesian};
            \addlegendentry{$p=3$}
            \addplot[fill=BarColorThree] table[x=label, y=p4] {\datatablecartesian};
            \addlegendentry{$p=4$}
            \addplot[fill=BarColorFour] table[x=label, y=p7] {\datatablecartesian};
            \addlegendentry{$p=7$}
        \end{axis}
    \end{tikzpicture}
    \caption{Cartesian mesh}
    \label{fig:throughput-cartesian-pmg}
\end{subfigure}
\hfill
\begin{subfigure}[b]{0.49\textwidth}

    \ReadAndProcessPatchFile{results/cartesian/pmg_p2.csv}{Two}

    \ReadAndProcessPatchFile{results/cartesian/pmg_p3.csv}{Three}

    \ReadAndProcessPatchFile{results/cartesian/pmg_p4.csv}{Four}
    \ReadAndProcessPatchFile{results/cartesian/pmg_p7.csv}{Seven}
    \centering
    \begin{tikzpicture}
        \pgfplotstableread[col sep=comma]{
        label,p2,p3,p4,p7
        {Smooth},\TwoPresmooth,\ThreePresmooth,\FourPresmooth,\SevenPresmooth
        {Restrict},\TwoRestrict,\ThreeRestrict,\FourRestrict,\SevenRestrict

        {Residual},\TwoResidual,\ThreeResidual,\FourResidual,\SevenResidual
        {Prolongate},\TwoProlongateThrough,\ThreeProlongateThrough,\FourProlongateThrough,\SevenProlongateThrough
        }\datatabledistorted

        \begin{axis}[
                width=\linewidth,
                height=6.2cm,
                ybar,
                bar width=5pt,
                enlarge x limits=0.25,
                symbolic x coords={{Smooth}, {Restrict}, {Residual}, {Prolongate}},
                xtick=data,
                xticklabel style={align=center, font=\small},
                ymin=0, ymax=20,
                ymajorgrids=true,
                grid style=dashed,
                legend style={at={(0.5,1.05)}, anchor=south, legend columns=-1, font=\small},
            ]
            \addplot[fill=BarColorOne] table[x=label, y=p2] {\datatabledistorted};
            \addlegendentry{$p=2$}
            \addplot[fill=BarColorTwo] table[x=label, y=p3] {\datatabledistorted};
            \addlegendentry{$p=3$}
            \addplot[fill=BarColorThree] table[x=label, y=p4] {\datatabledistorted};
            \addlegendentry{$p=4$}
            \addplot[fill=BarColorFour] table[x=label, y=p7] {\datatabledistorted};
            \addlegendentry{$p=7$}
        \end{axis}
    \end{tikzpicture}
    \caption{Distorted mesh}
    \label{fig:throughput-distorted-pmg}
\end{subfigure}
  \caption{Throughput of individual components of the p-multigrid application. The y-axis shows
    performance in GFLOP/s.}
  \label{fig:throughput-pmg}
\end{figure}

\paragraph{CPU thoughput}
To complement the timing analysis, we examine the computational throughput of the p-multigrid solver's internal
components, measured in GFLOP/s. High values indicate that the algorithm is compute-bound and effectively leveraging
the processor's vector units. Figure~\ref{fig:throughput-pmg} details this metric for the four main phases of the
V-cycle: smoothing, residual computation, restriction, and prolongation.

On the Cartesian mesh (Figure~\ref{fig:throughput-cartesian}), the residual computation stands out as the most
efficient kernel. Its throughput scales well with the polynomial degree, rising from approximately 8.8 GFLOP/s at $p=2$
to over 17.5 GFLOP/s at $p=7$. This confirms that the matrix-free operator evaluation, which dominates the residual
step, successfully saturates the floating-point units as the arithmetic intensity increases. In contrast, the
inter-grid transfer operators (restriction and prolongation) and the smoothing step exhibit significantly lower These
operations reach only about 4 GFLOP/s at $p=7$ and are clearly memory-bound. Unlike the operator evaluation, the
inter-grid transfer kernels perform simple streaming updates or interpolations and do not exploit SIMD. There is
therefore room for optimization: redesigning these transfers to be SIMD-friendly to be explicitly vectorized should
increase their throughput.

The results for the distorted mesh (Figure~\ref{fig:throughput-distorted}) reveal an important characteristic of the
implementation. Despite the increased mathematical complexity required to handle the deformed geometry—specifically the
full tensor contractions for the Jacobian—the throughput of the residual computation remains high, reaching
approximately 18.4 GFLOP/s at $p=7$. This is comparable to, and even slightly exceeds, the Cartesian case in terms of
raw floating-point operations per second. This indicates that while the distorted mesh requires more operations (and
thus takes longer in wall-clock time), the solver processes this additional workload efficiently. The implementation
remains compute-bound in this settings, successfully hiding the latency of the extra data fetches required for the
variable geometric factors behind a dense wall of arithmetic instructions.

\subsection{Comparison of Local Solver Strategies}
The choice of the local solver is a critical design decision for patch-based smoothers. As discussed in the
introduction, several strategies exist, ranging from simple iterative methods to sophisticated direct solvers. To
quantify the trade-offs involved, we compare the performance of our proposed p-multigrid variants against the widely
used fast diagonalization method~\cite{WitteArndtKanschat21}. Fast diagonalization is particularly relevant as a
baseline because it provides an efficient, direct inversion for tensor-product operators on Cartesian grids,
representing a "gold standard" for performance in ideal geometric conditions.

Table~\ref{tab:local-solver-strategies} presents a detailed comparison of these strategies for a single patch
application at a polynomial degree of $p=7$. We evaluate the standard V-cycle and the optimized half V-cycle, both with
the basic Jacobi smoother and the more advanced Cartesian-reinforced smoother. We report the wall-clock time relative
to the evaluation of the patch operator $A u$, providing a normalized metric that isolates algorithmic efficiency from
hardware specifics. Additionally, we track the total number of floating-point operations (FLOPs) and the computational
throughput to assess how well each method utilizes the hardware.

\begin{table}[htbp]
  \centering
  \footnotesize
  \caption{Comparison of different local solver strategies for a single cartesian patch application at polynomial
    degree $p=7$.
    The reported values are  normalized by the evaluation of the patch operator application $A u$ which was being
    evaluated at 22.7 GFLOP/s and required 491,530 floating-point operations.
  }
  \label{tab:local-solver-strategies}
  \begin{tabular}{lccc}
    \toprule
    \textbf{Local Solver Strategy}             & \textbf{Relative Time} & \textbf{Relative FLOPs} & \textbf{Throughput
      (GFLOP/s)}
    \\
    \midrule
    V-cycle (Jacobi)                           & 3.39                   & 2.27                    & 15.1
    \\
    Half V-cycle (Jacobi)                      & 2.01                   & 1.17                    & 13.2
    \\
    \addlinespace
    V-cycle (Cartesian-reinforced Jacobi)      & 8.68                   & 3.67                    & 9.59
    \\
    Half V-cycle (Cartesian-reinforced Jacobi) & 4.67                   & 1.87                    & 9.08
    \\
    \addlinespace
    \myrowcolor
    Fast Diagonalization                       & 2.60                   & 0.68                    & 5.95
    \\
    \bottomrule
  \end{tabular}
\end{table}

The results in Table~\ref{tab:local-solver-strategies} highlight the efficiency gains of the proposed optimization. The
standard V-cycle with a Jacobi smoother is already competitive, but the half V-cycle reduces the relative execution
time by approximately 40\% (from 3.39 to 2.01 times the cost of an operator evaluation), while maintaining a high
computational throughput of 13.2 GFLOP/s. Notably, the half V-cycle with Jacobi smoothing is faster than the fast
diagonalization method (2.01 vs. 2.60 relative time), despite performing more floating-point operations (1.17 vs. 0.68
relative FLOPs). This counter-intuitive result is explained by the significantly higher throughput of the matrix-free
p-multigrid operations compared to the memory-bound nature of the fast diagonalization on modern CPUs. While the
Cartesian-reinforced smoother offers better convergence properties, it comes at a higher computational cost due to its
lower throughput, making the simple Jacobi half V-cycle the most time-efficient choice for this configuration.

\subsection{Performance}
\label{sec:performance_parallel}

Finally, we transition from single-core micro-benchmarks on small problems to full-scale simulations running on a
multi-core system. The following results were obtained on a dual-socket node equipped with AMD EPYC 7282 processors,
providing a total of 32 physical cores. To handle the data dependencies inherent in the multiplicative patch smoother
in a parallel environment, we employ a graph coloring approach. The patches are colored using the DSATUR
algorithm~\cite{brelaz1979DSatur} such that no two patches sharing a degree of freedom share the same color. The
smoother is then applied color by color, with patches of the same color processed in parallel. While more advanced
scheduling algorithms exist to maximize cache locality~\cite{wichrowski2025smoothers}, the inherent data locality of
the patch-based method is already high. In this context, thread synchronization overhead becomes the primary concern,
and the simple coloring scheme provides a robust baseline for parallel execution.

The problem sizes for these benchmarks are detailed in Table~\ref{tab:dofs-patches}. They range from coarse meshes,
where synchronization overheads may be visible, to fine meshes with nearly 100 million degrees of freedom, where the
system is fully loaded.

\begin{table}[htbp]
  \centering
  \caption{Problem sizes for the multithreaded benchmarks. The number of patches is independent of the polynomial
    degree $p$.}
  \label{tab:dofs-patches}
  \begin{tabular}{lrrrrr}
    \toprule
    Refinements & Patches & \multicolumn{4}{c}{Number of DoFs}                                       \\
                &         & $p=2$                              & $p=3$     & $p=4$      & $p=7$      \\
    \midrule
    2           & 342     & 4,913                              & 15,625    & 35,937     & 185,193    \\
    \myrowcolor
    3           & 3,375   & 35,937                             & 117,649   & 274,625    & 1,442,897  \\
    4           & 29,791  & 274,625                            & 912,673   & 2,146,689  & 11,390,625 \\
    \myrowcolor
    5           & 250,047 & 2,146,689                          & 7,189,057 & 16,974,593 & 90,518,849 \\ \bottomrule
  \end{tabular}
\end{table}

\paragraph{Runtime Scaling}
We evaluate the performance of the parallel patch smoother relative to a global matrix-free operator evaluation ($Au$).
In a 3D hexahedral mesh, each cell belongs to eight vertex patches. Consequently, a naive lower bound for the cost of
one smoother application would be roughly eight times the cost of a global operator evaluation ($8 \times Au$), plus
the overhead of the local solves and data movement.

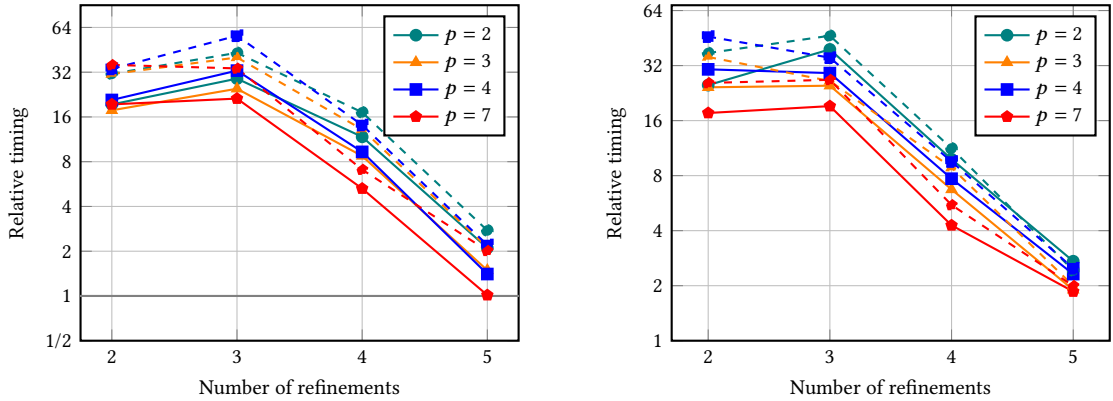
\begin{figure}[htbp]
  \centering
  \begin{tikzpicture}
    \begin{scope}
        \begin{axis}[
                paperplot,
                xlabel={Number of refinements},
                ylabel={Relative timing},
                xtick={0,1,2,3,4,5,6,7,8},
                ymode=log,
                ytick={0.5,1,2,4,8,16,32,64},
                yticklabels={1/2, 1,2,4,8,16,32,64},
                xmin=1.75, xmax=5.25,
                ymin=0.5, log basis x=2, legend pos=north east ] \pgfplotsset {cycle list set=0}

            \addplot+[gray, thick, forget plot] coordinates {(1,1) (6,1)};

            \addplot+[thick] table[ col sep=comma,
                    x=refinements,
                    y expr=\thisrow{smoother_avg_wall} / \thisrow{vmult_avg_wall},
                    restrict expr to domain={\thisrow{degree}}{1.9:2.1},
                    unbounded coords=discard ]
                {results/cartesian/threaded_n1.csv}; \addlegendentry{$p=2$}

            \addplot+[ thick] table[
                    col sep=comma,
                    x=refinements,
                    y expr=\thisrow{smoother_avg_wall} / \thisrow{vmult_avg_wall},
                    restrict expr to domain={\thisrow{degree}}{2.9:3.1},
                    unbounded coords=discard
                ] {results/cartesian/threaded_n1.csv};
            \addlegendentry{$p=3$}

            \addplot+[ thick] table[
                    col sep=comma,
                    x=refinements,
                    y expr=\thisrow{smoother_avg_wall} / \thisrow{vmult_avg_wall},
                    restrict expr to domain={\thisrow{degree}}{3.9:4.1},
                    unbounded coords=discard
                ] {results/cartesian/threaded_n1.csv};
            \addlegendentry{$p=4$}

            \addplot+[thick] table[
                    col sep=comma,
                    x=refinements,
                    y expr=\thisrow{smoother_avg_wall} / \thisrow{vmult_avg_wall},
                    restrict expr to domain={\thisrow{degree}}{6.9:7.1},
                    unbounded coords=discard
                ] {results/cartesian/threaded_n1.csv};
            \addlegendentry{$p=7$}

            \pgfplotsset {cycle list set=0}

            \addplot+[thick,dashed] table[ col sep=comma, x=refinements,
                    y expr=\thisrow{smoother_avg_wall} / \thisrow{vmult_avg_wall},
                    restrict expr to domain={\thisrow{degree}}{1.9:2.1},
                    unbounded coords=discard ]
                {results/cartesian/threaded_n2.csv};

            \addplot+[ thick,dashed] table[
                    col sep=comma,
                    x=refinements,
                    y expr=\thisrow{smoother_avg_wall} / \thisrow{vmult_avg_wall},
                    restrict expr to domain={\thisrow{degree}}{2.9:3.1},
                    unbounded coords=discard
                ] {results/cartesian/threaded_n2.csv};

            \addplot+[ thick,dashed] table[
                    col sep=comma,
                    x=refinements,
                    y expr=\thisrow{smoother_avg_wall} / \thisrow{vmult_avg_wall},
                    restrict expr to domain={\thisrow{degree}}{3.9:4.1},
                    unbounded coords=discard
                ] {results/cartesian/threaded_n2.csv};

            \addplot+[thick,dashed] table[
                    col sep=comma,
                    x=refinements,
                    y expr=\thisrow{smoother_avg_wall} / \thisrow{vmult_avg_wall},
                    restrict expr to domain={\thisrow{degree}}{6.9:7.1},
                    unbounded coords=discard
                ] {results/cartesian/threaded_n2.csv};
        \end{axis}
    \end{scope}

    \begin{scope}[xshift=0.52\columnwidth]
        \begin{axis}[
                paperplot,
                xlabel={Number of refinements},
                ylabel={Relative timing},
                xtick={0,1,2,3,4,5,6,7,8},
                ymode=log,
                ytick={1,2,4,8,16,32,64},
                yticklabels={1,2,4,8,16,32,64},
                ymin=1, log basis x=2,
                legend pos=north east ] \pgfplotsset {cycle list set=0}

            \addplot+[thick] table[ col sep=comma, x=refinements,
                    y expr=\thisrow{smoother_avg_wall} / \thisrow{vmult_avg_wall},
                    restrict expr to domain={\thisrow{degree}}{1.9:2.1},
                    unbounded coords=discard ]
                {results/affine/threaded_n1.csv}; \addlegendentry{$p=2$}

            \addplot+[ thick] table[
                    col sep=comma,
                    x=refinements,
                    y expr=\thisrow{smoother_avg_wall} / \thisrow{vmult_avg_wall},
                    restrict expr to domain={\thisrow{degree}}{2.9:3.1},
                    unbounded coords=discard
                ] {results/affine/threaded_n1.csv};
            \addlegendentry{$p=3$}

            \addplot+[ thick] table[
                    col sep=comma,
                    x=refinements,
                    y expr=\thisrow{smoother_avg_wall} / \thisrow{vmult_avg_wall},
                    restrict expr to domain={\thisrow{degree}}{3.9:4.1},
                    unbounded coords=discard
                ] {results/affine/threaded_n1.csv};
            \addlegendentry{$p=4$}

            \addplot+[thick] table[
                    col sep=comma,
                    x=refinements,
                    y expr=\thisrow{smoother_avg_wall} / \thisrow{vmult_avg_wall},
                    restrict expr to domain={\thisrow{degree}}{6.9:7.1},
                    unbounded coords=discard
                ] {results/affine/threaded_n1.csv};
            \addlegendentry{$p=7$}

            \pgfplotsset {cycle list set=0}

            \addplot+[thick,dashed] table[ col sep=comma, x=refinements,
                    y expr=\thisrow{smoother_avg_wall} / \thisrow{vmult_avg_wall},
                    restrict expr to domain={\thisrow{degree}}{1.9:2.1},
                    unbounded coords=discard ]
                {results/affine/threaded_n2.csv};

            \addplot+[ thick,dashed] table[
                    col sep=comma,
                    x=refinements,
                    y expr=\thisrow{smoother_avg_wall} / \thisrow{vmult_avg_wall},
                    restrict expr to domain={\thisrow{degree}}{2.9:3.1},
                    unbounded coords=discard
                ] {results/affine/threaded_n2.csv};

            \addplot+[ thick,dashed] table[
                    col sep=comma,
                    x=refinements,
                    y expr=\thisrow{smoother_avg_wall} / \thisrow{vmult_avg_wall},
                    restrict expr to domain={\thisrow{degree}}{3.9:4.1},
                    unbounded coords=discard
                ] {results/affine/threaded_n2.csv};

            \addplot+[thick,dashed] table[
                    col sep=comma,
                    x=refinements,
                    y expr=\thisrow{smoother_avg_wall} / \thisrow{vmult_avg_wall},
                    restrict expr to domain={\thisrow{degree}}{6.9:7.1},
                    unbounded coords=discard
                ] {results/affine/threaded_n2.csv};
        \end{axis}
    \end{scope}
\end{tikzpicture}
  \caption{
    Relative wall-clock time of the patch smoother application compared to a single global operator evaluation
    ($A u$) using a \emph{full V-cycle} local solver. Left: Cartesian mesh, right: distorted mesh. Solid lines: 1 V-cycle, dashed lines: 2 V-cycles. Runs were multithreaded on 32 cores.
  }
  \label{fig:threads-timing}
\end{figure}

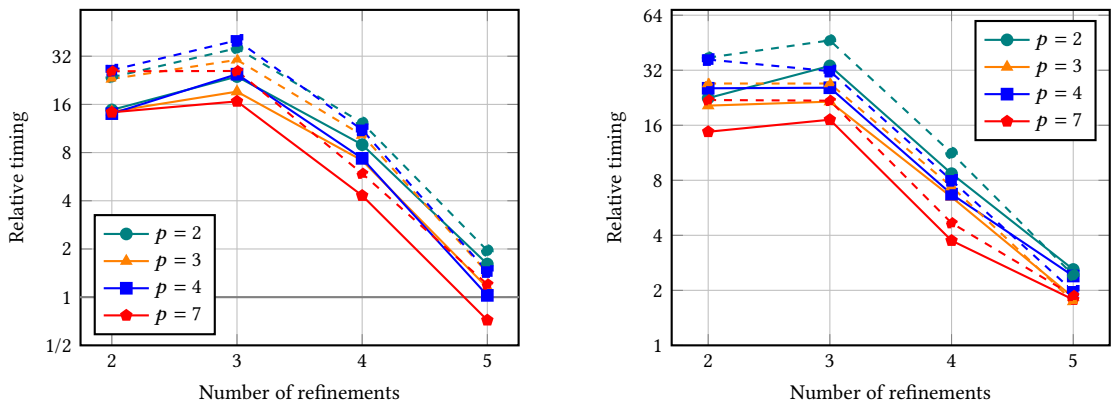
\begin{figure}[htbp]
  \centering
  \begin{tikzpicture}
    \begin{scope}
        \begin{axis}[
                paperplot,
                xlabel={Number of refinements},
                ylabel={Relative timing},
                xtick={0,1,2,3,4,5,6,7,8},
                ymode=log,
                ytick={1/2, 1,2,4,8,16,32,64},
                yticklabels={1/2, 1,2,4,8,16,32,64},
                xmin=1.75, xmax=5.25,
                ymin=1/2, log basis x=2, legend pos=south west ] \pgfplotsset {cycle list set=0}

            \addplot+[gray, thick, forget plot] coordinates {(1,1) (6,1)};

            \addplot+[thick] table[ col sep=comma,
                    x=refinements,
                    y expr=\thisrow{smoother_avg_wall} / \thisrow{vmult_avg_wall},
                    restrict expr to domain={\thisrow{degree}}{1.9:2.1},
                    unbounded coords=discard ]
                {results/cartesian/threaded_half_n1.csv}; \addlegendentry{$p=2$}

            \addplot+[ thick] table[
                    col sep=comma,
                    x=refinements,
                    y expr=\thisrow{smoother_avg_wall} / \thisrow{vmult_avg_wall},
                    restrict expr to domain={\thisrow{degree}}{2.9:3.1},
                    unbounded coords=discard
                ] {results/cartesian/threaded_half_n1.csv};
            \addlegendentry{$p=3$}

            \addplot+[ thick] table[
                    col sep=comma,
                    x=refinements,
                    y expr=\thisrow{smoother_avg_wall} / \thisrow{vmult_avg_wall},
                    restrict expr to domain={\thisrow{degree}}{3.9:4.1},
                    unbounded coords=discard
                ] {results/cartesian/threaded_half_n1.csv};
            \addlegendentry{$p=4$}

            \addplot+[thick] table[
                    col sep=comma,
                    x=refinements,
                    y expr=\thisrow{smoother_avg_wall} / \thisrow{vmult_avg_wall},
                    restrict expr to domain={\thisrow{degree}}{6.9:7.1},
                    unbounded coords=discard
                ] {results/cartesian/threaded_half_n1.csv};
            \addlegendentry{$p=7$}

            \pgfplotsset {cycle list set=0}

            \addplot+[thick,dashed] table[ col sep=comma, x=refinements,
                    y expr=\thisrow{smoother_avg_wall} / \thisrow{vmult_avg_wall},
                    restrict expr to domain={\thisrow{degree}}{1.9:2.1},
                    unbounded coords=discard ]
                {results/cartesian/threaded_half_n2.csv};

            \addplot+[ thick,dashed] table[
                    col sep=comma,
                    x=refinements,
                    y expr=\thisrow{smoother_avg_wall} / \thisrow{vmult_avg_wall},
                    restrict expr to domain={\thisrow{degree}}{2.9:3.1},
                    unbounded coords=discard
                ] {results/cartesian/threaded_half_n2.csv};

            \addplot+[ thick,dashed] table[
                    col sep=comma,
                    x=refinements,
                    y expr=\thisrow{smoother_avg_wall} / \thisrow{vmult_avg_wall},
                    restrict expr to domain={\thisrow{degree}}{3.9:4.1},
                    unbounded coords=discard
                ] {results/cartesian/threaded_half_n2.csv};

            \addplot+[thick,dashed] table[
                    col sep=comma,
                    x=refinements,
                    y expr=\thisrow{smoother_avg_wall} / \thisrow{vmult_avg_wall},
                    restrict expr to domain={\thisrow{degree}}{6.9:7.1},
                    unbounded coords=discard
                ] {results/cartesian/threaded_half_n2.csv};
        \end{axis}
    \end{scope}

    \begin{scope}[xshift=0.52\columnwidth]
        \begin{axis}[
                paperplot,
                xlabel={Number of refinements},
                ylabel={Relative timing},
                xtick={0,1,2,3,4,5,6,7,8},
                ymode=log,
                ytick={1,2,4,8,16,32,64},
                yticklabels={1,2,4,8,16,32,64},
                ymin=1, log basis x=2,
                legend pos=north east ] \pgfplotsset {cycle list set=0}

            \addplot+[thick] table[ col sep=comma, x=refinements,
                    y expr=\thisrow{smoother_avg_wall} / \thisrow{vmult_avg_wall},
                    restrict expr to domain={\thisrow{degree}}{1.9:2.1},
                    unbounded coords=discard ]
                {results/affine/threaded_half_n1.csv}; \addlegendentry{$p=2$}

            \addplot+[ thick] table[
                    col sep=comma,
                    x=refinements,
                    y expr=\thisrow{smoother_avg_wall} / \thisrow{vmult_avg_wall},
                    restrict expr to domain={\thisrow{degree}}{2.9:3.1},
                    unbounded coords=discard
                ] {results/affine/threaded_half_n1.csv};
            \addlegendentry{$p=3$}

            \addplot+[ thick] table[
                    col sep=comma,
                    x=refinements,
                    y expr=\thisrow{smoother_avg_wall} / \thisrow{vmult_avg_wall},
                    restrict expr to domain={\thisrow{degree}}{3.9:4.1},
                    unbounded coords=discard
                ] {results/affine/threaded_half_n1.csv};
            \addlegendentry{$p=4$}

            \addplot+[thick] table[
                    col sep=comma,
                    x=refinements,
                    y expr=\thisrow{smoother_avg_wall} / \thisrow{vmult_avg_wall},
                    restrict expr to domain={\thisrow{degree}}{6.9:7.1},
                    unbounded coords=discard
                ] {results/affine/threaded_half_n1.csv};
            \addlegendentry{$p=7$}

            \pgfplotsset {cycle list set=0}

            \addplot+[thick,dashed] table[ col sep=comma, x=refinements,
                    y expr=\thisrow{smoother_avg_wall} / \thisrow{vmult_avg_wall},
                    restrict expr to domain={\thisrow{degree}}{1.9:2.1},
                    unbounded coords=discard ]
                {results/affine/threaded_n2.csv};

            \addplot+[ thick,dashed] table[
                    col sep=comma,
                    x=refinements,
                    y expr=\thisrow{smoother_avg_wall} / \thisrow{vmult_avg_wall},
                    restrict expr to domain={\thisrow{degree}}{2.9:3.1},
                    unbounded coords=discard
                ] {results/affine/threaded_half_n2.csv};

            \addplot+[ thick,dashed] table[
                    col sep=comma,
                    x=refinements,
                    y expr=\thisrow{smoother_avg_wall} / \thisrow{vmult_avg_wall},
                    restrict expr to domain={\thisrow{degree}}{3.9:4.1},
                    unbounded coords=discard
                ] {results/affine/threaded_half_n2.csv};

            \addplot+[thick,dashed] table[
                    col sep=comma,
                    x=refinements,
                    y expr=\thisrow{smoother_avg_wall} / \thisrow{vmult_avg_wall},
                    restrict expr to domain={\thisrow{degree}}{6.9:7.1},
                    unbounded coords=discard
                ] {results/affine/threaded_half_n2.csv};
        \end{axis}
    \end{scope}
\end{tikzpicture}
  \caption{
    Relative wall-clock time of the patch smoother application compared to a single global operator evaluation
    ($A u$) using the optimized \emph{Half V-cycle} local solver. Left: Cartesian mesh, right: distorted mesh. Solid lines: 1 half V-cycle, dashed lines: 2 half V-cycles. Runs were multithreaded on 32 cores.
  }
  \label{fig:half-vcycle-threads-timing}
\end{figure}

Figure~\ref{fig:threads-timing} presents the relative wall-clock time for the standard \emph{full V-cycle} local
solver, while Figure~\ref{fig:half-vcycle-threads-timing} presents the results for the optimized \emph{half V-cycle}.
For coarser meshes, the relative timings are significantly higher than the baseline factor of 8. This is attributed to
the limited parallel work available per color and the relative impact of thread synchronization barriers between
colors.

However, as the mesh is refined, the relative cost decreases dramatically. This phenomenon is most pronounced for the
Cartesian mesh at $p=7$. On the finest level, the cost of the entire smoother application drops below the theoretical
baseline of $8 \times Au$. This occurs because the global operator evaluation $Au$ becomes increasingly memory-bound on
large problems, with performance limited by the speed of streaming vectors from main memory. In contrast, the patch
smoother, while performing significantly more arithmetic operations, operates on smaller, cache-resident data chunks
per patch. This allows it to achieve higher aggregate arithmetic intensity and better utilize the available memory
bandwidth.

Noteworthy, for $p=5$ the application of smoother becomes less expensive than evaluation of $Au$. This indicates that
the current multithreaded implementation of the global matrix-free operator loop in \texttt{deal.II} may not be fully
optimized for this regime, potentially suffering from synchronization overheads or suboptimal bandwidth utilization
that our explicit patch-based scheduling manages to avoid.

Comparing Figure~\ref{fig:threads-timing} and Figure~\ref{fig:half-vcycle-threads-timing} reveals a modest performance
improvement from the half V-cycle optimization. Because the smoother's execution time is partially dominated by the
fixed overheads of data movement and setup, the reduction in local arithmetic intensity does not translate linearly to
the total runtime.

\begin{figure}[htbp]
  \centering
  \begin{tikzpicture}
    \begin{scope}
        \begin{axis}[
                paperplot,
                xlabel={Number of refinements},
                ylabel={Memory [doubles/DoF]},
                xtick={0,1,2,3,4,5,6,7,8},
                ymode=log,
                ytick={1,2,4,8,16,32,64},
                yticklabels={1,2,4,8,16,32,64},
                ymin=2, log basis x=2, legend pos=south east,  legend pos=south east, legend columns=2 ] \pgfplotsset {cycle list set=0}

            \addplot+[thick] table[ col sep=comma,
                    x=refinements,
                    y expr=\thisrow{mem_bytes} / 8 / \thisrow{dofs},
                    restrict expr to domain={\thisrow{degree}}{1.9:2.1},
                    unbounded coords=discard ]
                {results/cartesian/threaded_half_n1.csv}; \addlegendentry{$p=2$}

            \addplot+[ thick] table[
                    col sep=comma,
                    x=refinements,
                    y expr=\thisrow{mem_bytes} / 8 / \thisrow{dofs},
                    restrict expr to domain={\thisrow{degree}}{2.9:3.1},
                    unbounded coords=discard
                ] {results/cartesian/threaded_half_n1.csv};
            \addlegendentry{$p=3$}

            \addplot+[ thick] table[
                    col sep=comma,
                    x=refinements,
                    y expr=\thisrow{mem_bytes} / 8 / \thisrow{dofs},
                    restrict expr to domain={\thisrow{degree}}{3.9:4.1},
                    unbounded coords=discard
                ] {results/cartesian/threaded_half_n1.csv};
            \addlegendentry{$p=4$}

            \addplot+[thick] table[
                    col sep=comma,
                    x=refinements,
                    y expr=\thisrow{mem_bytes} / 8 / \thisrow{dofs},
                    restrict expr to domain={\thisrow{degree}}{6.9:7.1},
                    unbounded coords=discard
                ] {results/cartesian/threaded_half_n1.csv};
            \addlegendentry{$p=7$}

        \end{axis}
    \end{scope}

    \begin{scope}[xshift=0.52\columnwidth]
        \begin{axis}[
                paperplot,
                xlabel={Number of refinements},
                ylabel={Memory [doubles/DoF]},
                xtick={0,1,2,3,4,5,6,7,8},
                ymode=log,
                ytick={1,2,4,8,16,32,64},
                yticklabels={1,2,4,8,16,32,64},
                ymin=16, log basis x=2,
                legend pos=south east, legend columns=2 ] \pgfplotsset {cycle list set=0}

            \addplot+[thick] table[ col sep=comma, x=refinements,
                    y expr=\thisrow{mem_bytes} / 8 / \thisrow{dofs},
                    restrict expr to domain={\thisrow{degree}}{1.9:2.1},
                    unbounded coords=discard ]
                {results/affine/threaded_half_n1.csv}; \addlegendentry{$p=2$}

            \addplot+[ thick] table[
                    col sep=comma,
                    x=refinements,
                    y expr=\thisrow{mem_bytes} / 8 / \thisrow{dofs},
                    restrict expr to domain={\thisrow{degree}}{2.9:3.1},
                    unbounded coords=discard
                ] {results/affine/threaded_half_n1.csv};
            \addlegendentry{$p=3$}

            \addplot+[ thick] table[
                    col sep=comma,
                    x=refinements,
                    y expr=\thisrow{mem_bytes} / 8 / \thisrow{dofs},
                    restrict expr to domain={\thisrow{degree}}{3.9:4.1},
                    unbounded coords=discard
                ] {results/affine/threaded_half_n1.csv};
            \addlegendentry{$p=4$}

            \addplot+[thick] table[
                    col sep=comma,
                    x=refinements,
                    y expr=\thisrow{mem_bytes} / 8 / \thisrow{dofs},
                    restrict expr to domain={\thisrow{degree}}{6.9:7.1},
                    unbounded coords=discard
                ] {results/affine/threaded_half_n1.csv};
            \addlegendentry{$p=7$}

        \end{axis}
    \end{scope}
\end{tikzpicture}
  \caption{
    Memory footprint of the patch smoother application (in doubles per DoF) for varying problem sizes and polynomial degrees $p$. Left: Cartesian mesh, right: distorted mesh.
  }
  \label{fig:memory-footprint}
\end{figure}
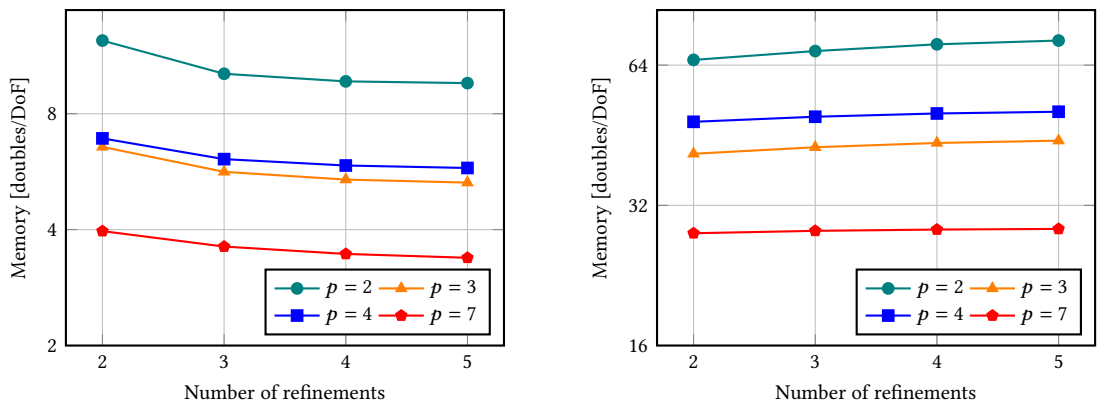

\paragraph{Memory Footprint}
The memory footprint of the patch smoother is a critical constraint for high-order methods. Our implementation requires
storing the inverse diagonal of the operator on each p-level for the Jacobi smoother. Figure~\ref{fig:memory-footprint}
shows the memory consumption measured in doubles per degree of freedom.

The results confirm the $\mathcal{O}(p^d)$ scaling of the storage requirements. Notably, the memory cost per DoF
\emph{decreases} or remains constant as $p$ increases. For the Cartesian case, the footprint is minimal (approx. 2-4
doubles/DoF), as the operator evaluation does not require storing geometric factors. For the distorted mesh, the
footprint is higher (approx. 16-32 doubles/DoF) because the matrix-free operator evaluation within the smoother
requires access to the stored Jacobian determinants and inverse Jacobians at quadrature points. However, importantly,
the memory usage per DoF remains stable as the number of refinements increases, confirming that the method scales
perfectly to large problems without memory overhead blow-up.

\paragraph{Measurement Note}
While the relative performance trends strongly suggest that the global operator evaluation becomes memory-bound on fine
meshes, we were unable to directly verify this hypothesis using hardware performance counters in the multi-threaded
benchmarks. The performance analysis tool LIKWID, which we successfully employed for the single-core micro-benchmarks,
relies on specific thread pinning that conflicted with the dynamic task scheduling of Intel's Threading Building Blocks
(TBB) used in our parallel implementation. Nevertheless, the observed scaling behavior is consistent with the
transition from a compute-bound to a memory-bound regime for the baseline operator, as analyzed
in~\cite{wichrowski2025smoothers}.

\section{Conclusion}

In this work, we have presented a blueprint for a matrix-free, p-multigrid patch smoother that resolves the classic
dilemma between rapid convergence and rapid execution. While vertex-patch smoothers are prized for their low iteration
counts, they have traditionally been penalized by a high computational cost per step. By employing a local p-multigrid
solver --- efficiently realized through sum-factorization and explicit SIMD vectorization --- we have successfully
reduced this cost to a level comparable with a simple Jacobi smoother for $p=7$ on the finest mesh, effectively
democratizing access to Schwarz-type robustness for high-order applications.

Our analysis of the "half V-cycle" optimization offers a practical lesson in modern hardware realities: while
eliminating the post-smoothing step significantly reduces the arithmetic operation count, the resulting speedup is
modest. This serves as a gentle reminder that in the regime of high-performance computing, saving FLOPs is only half
the battle; the other half is fighting the limitation of memory bandwidth. Nevertheless, the method maintains an
optimal $\mathcal{O}(p^d)$ memory footprint and scales favorably on fine grids, where the compute-intensive patch
smoother eventually outpaces the memory-bound global operator evaluation.

Ultimately, numerical analysts are often resigned to the trilemma of designing solvers that are \emph{Robust, Fast, and
  General} --- usually, one is forced to pick only two. While we cannot claim to have completely disproven the "No Free
Lunch" theorem, our results suggest that with enough matrix-free optimization, we have at least negotiated a very
significant discount on the bill.

\section*{Acknowledgments}

The generative AI (Gemini, ChatGPT, Claude) was used to assist in drafting and proofreading parts of this manuscript.
Any remaining errors are the author's responsibility.

\noindent The author warmly thanks his mother, Grażyna, and his brother, Wiktor, for their support. Special thanks to Micro, the family cat, for her comforting companionship during the writing of this paper.

\bibliographystyle{siamplain}
\bibliography{literature}

\end{document}